\documentclass{article}
\usepackage[latin9]{inputenc}
\usepackage{mathtools}
\usepackage{amsmath}
\usepackage{amssymb}

\makeatletter
\usepackage{amsthm}\usepackage{enumerate}
\usepackage{comment}
\numberwithin{equation}{section}
\newtheorem{defn}{Definition}[section]
\newtheorem{lem}[defn]{Lemma}\newtheorem{thm}[defn]{Theorem}\newtheorem{prop}[defn]{Proposition}\newtheorem{cor}[defn]{Corollary}\newtheorem{rem}[defn]{Remark}\newtheorem{rem*}[defn]{Remark}\textwidth=17cm
\makeatother

\begin{document}

\title{Extended Hamilton-Jacobi Theory, symmetries and integrability by
quadratures}

\author{S. Grillo, J.C. Marrero \& E. Padrón}
\maketitle
\begin{abstract}
In this paper, we study the extended Hamilton-Jacobi Theory in the
context of dynamical systems with symmetries. Given an action of a
Lie group $G$ on a manifold $M$ and a $G$-invariant vector field
$X$ on $M$, we construct complete solutions of the Hamilton-Jacobi
equation (HJE) related to $X$ (and a given fibration on $M$). We
do that along each open subset $U\subseteq M$ such that $\pi\left(U\right)$
has a manifold structure and $\pi_{\left|U\right.}:U\rightarrow\pi\left(U\right)$,
the restriction to $U$ of the canonical projection $\pi:M\rightarrow M/G$,
is a surjective submersion. If $X_{\left|U\right.}$ is not vertical
with respect to $\pi_{\left|U\right.}$, we show that such complete
solutions solve the \textit{reconstruction equations} related to $X_{\left|U\right.}$
and $G$, i.e., the equations that enable us to write the integral
curves of $X_{\left|U\right.}$ in terms of those of its projection
on $\pi\left(U\right)$. On the other hand, if $X_{\left|U\right.}$
is vertical, we show that such complete solutions can be used to construct
(around some points of $U$) the integral curves of $X_{\left|U\right.}$
up to quadratures. To do that we give, for some elements $\xi$ of
the Lie algebra $\mathfrak{g}$ of $G$, an explicit expression up
to quadratures of the exponential curve $\exp\left(\xi\,t\right)$,
different to that appearing in the literature for matrix Lie groups.
In the case of compact and of semisimple Lie groups, we show that
such expression of $\exp\left(\xi\,t\right)$ is valid for all $\xi$
inside an open dense subset of $\mathfrak{g}$. 
\end{abstract}
\tableofcontents{}

\section{Introduction}

In the last few years, several generalizations of the \textit{classical}
Hamilton-Jacobi equation (HJE) have been developed for Hamiltonian
systems on different contexts: on symplectic, cosymplectic, contact,
Poisson and almost-Poisson manifolds, and also on Lie algebroids.
The resulting Hamilton-Jacobi theories were applied to nonholonomic
systems, dissipative and time-dependent Hamiltonian systems, reduced
systems by symmetries and Hamiltonian systems with external forces
\cite{bmmp,pepin-holo,pepin-noholo,lmm,del,hjp,ds,fer}. In all of
them, the following ingredients are present: $\left(1\right)$ a fibration
$\Pi\kern-2pt :\kern-2pt M\rightarrow N$ (i.e. a surjective submersion)
defined on the phase space $M$ of each system; $\left(2\right)$
the solutions of the generalized HJE, which we shall call $\Pi$-HJE,
given by sections $\sigma\kern-2pt :\kern-2pt N\rightarrow M$ of
such a fibration $\Pi$; $\left(3\right)$ the complete solutions
$\Sigma\kern-2pt :\kern-2pt N\times\Lambda\kern-1pt \rightarrow\kern-1pt M$,
given by local diffeomorphisms such that, for each $\lambda\in\Lambda$,
$\sigma_{\lambda}:=\Sigma(\cdot,\lambda)$ is a solution of the $\Pi$-HJE.
This clearly generalizes the classical situation \cite{am,gold},
where the involved fibration is the cotangent projection $\pi_{Q}:T^{\ast}Q\rightarrow Q$
of a manifold $Q$ and the solutions $\sigma\kern-2pt :\kern-2pt Q\rightarrow T^{\ast}Q$
are exact $1$-forms on $Q$.

In Ref. \cite{gp}, an extension to general (i.e. not necessarily
Hamiltonian) dynamical systems, of the previously mentioned Hamilton-Jacobi
theories, was carried out, focusing on the connection between complete
solutions and the integrability by quadratures of the involved systems.

The main aim of the present paper is to further study such an extended
theory in the context of dynamical systems with symmetry. Concretely,
given a general action $\rho:G\times M\rightarrow M$ (not necessarily
free or proper) of a Lie group $G$ on a manifold $M$ and a $G$-invariant
vector field $X$ on $M$ (with respect to $\rho$), we investigate
how to use $\rho$ to construct (local) fibrations $\Pi$ of $M$
and related solutions of the $\Pi$-HJE for $X$. We first show that,
around almost every point of $M$ (depending on the isotropy subgroups
of $G$), there exists a neighborhood $U$ such that the canonical
projection $\pi:M\rightarrow M/G$ restricted to $U$, namely $\pi_{\left|U\right.}:U\rightarrow\pi\left(U\right)$,
defines a fibration (even though $\rho$ is neither free nor proper).
Then we consider two kinds of vector fields: $\left(a\right)$ those
for which $X_{\left|U\right.}$ is not vertical with respect to $\pi_{\left|U\right.}$,
which we call \textit{horizontal}, and $\left(b\right)$ the vertical
ones. For the horizontal vector fields we show that, related to the
action $\rho$, there exists a submersion $\Theta$ transverse to
$\pi_{\left|U\right.}$ (which plays the role of a flat principal
connection) such that 
\begin{equation}
\Sigma\coloneqq\left(\pi_{\left|U\right.},\Theta\right)^{-1}:\pi\left(U\right)\times\Lambda\rightarrow U,\label{csh}
\end{equation}
with $\Lambda$ a submanifold of $G$, is a complete solution of the
$\pi_{\left|U\right.}$-HJE for $X_{\left|U\right.}$. Such a $\Sigma$
can be seen as a solution of a \textit{reconstruction problem}, in
the sense that, if we know the integral curves $\gamma\left(t\right)$
of the projected vector field $Y$ of $X_{\left|U\right.}$ on $\pi\left(U\right)$,
then the integral curves of $X_{\left|U\right.}$ are given by $\Gamma\left(t\right)=\Sigma\left(\gamma\left(t\right),\lambda\right)$,
with $\lambda\in\Lambda$. For the vertical vector fields, we show
that we can construct up to quadratures a submersion $\Theta$ transverse
to $\pi_{\left|U\right.}$ such that 
\begin{equation}
\Sigma\coloneqq\left(\Theta,\pi_{\left|U\right.}\right)^{-1}:N\times\pi\left(U\right)\rightarrow U,\label{csv}
\end{equation}
with $N$ a submanifold of $G$, is a complete solution of the $\Theta$-HJE
for $X_{\left|U\right.}$. Moreover, we prove that the integral curves
of $X_{\left|U\right.}$ also can be constructed up to quadratures
around some points of $M$. To do that, we first show that the exponential
curves $t\mapsto\exp\left(\xi\,t\right)$ of $G$, for some elements
$\xi$ of its Lie algebra $\mathfrak{g}$, can be constructed up to
quadratures. As it is well-known, there exists several explicit expressions
of $\exp\left(\xi\,t\right)$, unless for matrix Lie groups. What
we are giving here is an alternative expression for such curves, valid
also for non-matrix Lie groups. In the case in which $G$ is semisimple
or compact, we show that such an expression is valid for all $\xi$
in an open dense subset of $\mathfrak{g}$.

\bigskip{}

The paper is organized as follows. In Section $2$, we make a brief
review of the extended Hamilton-Jacobi Theory appearing in \cite{gp,gp2}.
We also present a result that ensures, in the presence of a complete
solution and in the context of symplectic manifolds, the integrability
by quadratures of general vector fields. It extends a result proved
in \cite{gp} for Hamiltonian vector fields only. In Section $3$,
given a dynamical system with symmetry, we construct the complete
solutions \eqref{csh} and \eqref{csv} for horizontal and vertical
vector fields, respectively. In Section $4$, we show the intimate
relationship that there exists between the complete solutions of a
horizontal (and invariant) vector field and the so-called reconstruction
processes. Finally, in Section $5$, using the above mentioned result
of Section $2$, we show that the exponential curves $t\mapsto\exp\left(\xi\,t\right)$
of $G$, for some points $\xi\in\mathfrak{g}$, can be constructed
up to quadratures. Then, based on that, we state sufficient conditions
under which a vertical (and invariant) vector field can also be integrable
up to quadratures.\bigskip{}

We assume that the reader is familiar with the main concepts of Differential
Geometry (see \cite{boot,kn,mrgm}) and with the basic ideas related
to Hamiltonian systems, Symplectic Geometry and Poisson Geometry in
the context of Geometric Mechanics (see for instance \cite{am,ar,lib,mr}).
We shall work in the smooth (i.e. $C^{\infty}$) category, focusing
exclusively on finite-dimensional smooth manifolds. Regarding the
terminology associated to the concept of ``integrability by quadratures,''
we shall adopt the following convention. We shall say that ``a function
$F:P\rightarrow Q$ \textbf{can be constructed up to quadratures},''
or simply ``\textbf{can be constructed},'' if its domain $P$ and
its values $F\left(p\right)$ (for all $p\in P$): 
\begin{itemize}
\item are simply known; 
\item they can be determined by making a finite number of arithmetic operations
(as the calculation of a determinant) and/or solving a finite set
of linear equations (which actually can be reduced to arithmetic operations); 
\item or they can be expressed in terms of the derivatives, primitives (i.e.
quadratures) and/or lateral inverses (by using the Implicit or Inverse
Function Theorem) of another known functions. 
\end{itemize}
When the function $F$ above is an integral curve $\Gamma$ of a vector
field and such a curve can be constructed up to quadratures, we shall
say that $\Gamma$ can be \textbf{integrated up to quadratures}, or
\textbf{by quadratures}.

\section{Preliminaries: complete solutions, first integrals and integrability}

\subsection{The extended Hamilton-Jacobi equation}

\label{comp}

Consider a manifold $M$, a vector field $X\in\mathfrak{X}\left(M\right)$
and a surjective submersion $\Pi:M\rightarrow N$ (\emph{ipso facto}
an open map). Related to this data (see \cite{gp}) we have the $\Pi$-\textbf{Hamilton-Jacobi
equation }($\Pi$-\textbf{HJE}) for $X$: 
\begin{equation}
\sigma_{\ast}\circ\Pi_{\ast}\circ X\circ\sigma=X\circ\sigma,\label{hjer}
\end{equation}
whose unknown is a section $\sigma:N\rightarrow M$ of $\Pi$ (\emph{ipso
facto} a closed map). If $\sigma$ solves the equation above, we shall
say that $\sigma$ is a \textbf{(global) solution of the }$\Pi$-\textbf{HJE
for} $X$. On the other hand, given an open subset $U\subseteq M$,
we shall call \textbf{local solution of} \textbf{the }$\Pi$-\textbf{HJE
for} $X$ \textbf{along} $U$ to any solution of the $\Pi_{\left|U\right.}$-HJE
for $X_{\left|U\right.}$. (Here, we are seeing $\Pi_{\left|U\right.}$
as a submersion onto $\Pi\left(U\right)$ and $X_{\left|U\right.}$
as an element of $\mathfrak{X}\left(U\right)$). Note that $\sigma$
is a solution of the $\Pi$-HJE for $X$ if and only if 
\begin{equation}
\sigma_{*}\circ X^{\sigma}=X\circ\sigma,\label{hjrel}
\end{equation}
where 
\begin{equation}
X^{\sigma}\coloneqq\Pi_{*}\circ X\circ\sigma,\label{Xs}
\end{equation}
i.e. the vector fields $X\in\mathfrak{X}\left(M\right)$ and $X^{\sigma}\in\mathfrak{X}\left(N\right)$
are $\sigma$-related. (Moreover, it can be shown that $\sigma$ is
a solution of \eqref{hjer} if and only if its image is an $X$-invariant
closed submanifold). This means that, in order to find the trajectories
of $X$ along the image of $\sigma$, we can look for the integral
curves of $X^{\sigma}$.

Given another manifold $\Lambda$ such that $\dim\Lambda+\dim N=\dim M$,
a \textbf{complete solution of the }$\Pi$-\textbf{HJE for }$X$ is
a surjective local diffeomorphism $\Sigma:N\times\Lambda\rightarrow M$
such that, for all $\lambda\in\Lambda$, 
\begin{equation}
\sigma_{\lambda}:=\Sigma\left(\cdot,\lambda\right):p\in N\longmapsto\Sigma\left(p,\lambda\right)\in M\label{psr}
\end{equation}
is a solution of the $\Pi$-HJE for $X$. The local version is obtained
by replacing $M$, $X$, $\Pi$ and $N$ by $U$, $X_{\left|U\right.}$,
$\Pi_{\left|U\right.}$ and $\Pi\left(U\right)$, respectively, being
$U$ an open subset of $M$. Each section $\sigma_{\lambda}$ is called
a \textbf{partial solution}. We showed in \cite{gp} that a (local)
complete solution $\Sigma$ exists around every point $m\in M$ for
which $X\left(m\right)\notin\mathsf{Ker}\Pi_{\ast,m}$.

Denoting by $\mathfrak{p}_{N}$ and $\mathfrak{p}_{\Lambda}$ the
projections of $N\times\Lambda$ onto $N$ and $\Lambda$, respectively,
it is easy to prove that a surjective local diffeomorphism $\Sigma$
is a complete solution if and only if 
\begin{equation}
\Pi\circ\Sigma=\mathfrak{p}_{N}\ \ \ \text{and}\ \ \ \Sigma_{\ast}\circ X^{\Sigma}=X\circ\Sigma,\label{Srelr}
\end{equation}
being $X^{\Sigma}\in\mathfrak{X}\left(N\times\Lambda\right)$ the
unique vector field on $N\times\Lambda$ satisfying 
\begin{equation}
\left(\mathfrak{p}_{N}\right)_{\ast}\circ X^{\Sigma}=\Pi_{\ast}\circ X\circ\Sigma\ \ \ \text{and}\ \ \ \left(\mathfrak{p}_{\Lambda}\right)_{\ast}\circ X^{\Sigma}=0.\label{Srelrr}
\end{equation}
Note that $X^{\Sigma}\left(p,\lambda\right)=\left(X^{\sigma_{\lambda}}\left(p\right),0\right)$,
with $X^{\sigma_{\lambda}}:=\Pi_{\ast}\circ X\circ\sigma_{\lambda}\in\mathfrak{X}\left(N\right)$,
so, in particular, 
\begin{equation}
\mathsf{Im}X^{\Sigma}\subseteq TN\times\left\{ 0\right\} .\label{xsin}
\end{equation}
Also, the fields $X$ and $X^{\Sigma}$ are $\Sigma$-related. This
implies that all the trajectories of $X$ can be obtained from those
of $X^{\Sigma}$. More precisely, since each integral curve of $X^{\Sigma}$
is clearly of the form $t\mapsto\left(\gamma\left(t\right),\lambda\right)\in N\times\Lambda$,
for some fixed point $\lambda\in\Lambda$ (see \eqref{xsin}), those
of $X$ are given by 
\begin{equation}
t\mapsto\Sigma\left(\gamma\left(t\right),\lambda\right)=\sigma_{\lambda}\left(\gamma\left(t\right)\right).\label{G}
\end{equation}
So, for each $\lambda$, we only need to find the curves $\gamma$,
which are the integral curves of the vector field $X^{\sigma_{\lambda}}\in\mathfrak{X}\left(N\right)$.

\subsection{The ``complete solutions - first integrals'' \textit{duality}}

\label{dual} Consider again a manifold $M$, a vector field $X\in\mathfrak{X}\left(M\right)$
and a surjective submersion $\Pi:M\rightarrow N$. We shall say that
a submersion $F:M\rightarrow\Lambda$ is a \textbf{first integrals
submersion} if 
\begin{equation}
\mathsf{Im}X\subseteq\mathsf{Ker}F_{*}.\label{fi}
\end{equation}

\begin{rem*} Note that, if $\Lambda=\mathbb{R}^{l}$, the components
$f_{1},...,f_{l}:M\rightarrow\mathbb{R}$ of $F$ define a set of
$l$ (functionally) independent first integrals, in the usual sense.
\end{rem*} Also, we shall say that $F$ is \textbf{transverse to}
$\Pi$ if 
\begin{equation}
TM=\mathsf{Ker}\Pi_{*}\oplus\mathsf{Ker}F_{*}.\label{tr}
\end{equation}
It was shown in \cite{gp} that, given a complete solution $\Sigma:N\times\Lambda\rightarrow M$
of the $\Pi$-HJE for $X$, we can construct around every point of
$M$ a neighborhood $U$ and a submersion $F:U\rightarrow\Lambda$
such that 
\begin{itemize}
\item $\mathsf{Im}X_{\left|U\right.}\subseteq\mathsf{Ker}F_{*}$ (first
integrals), 
\item $TU=\mathsf{Ker}\left(\Pi_{\left|U\right.}\right)_{*}\oplus\mathsf{Ker}F_{*}$
(transversality). 
\end{itemize}
In other words, from $\Sigma$ we have, around every point of $M$,
a first integrals submersion transverse to $\Pi$. The subset $U$
and the function $F$ are given by the formulae 
\begin{equation}
U\coloneqq\Sigma\left(V\right)\;\;\;\textrm{and \;\;\;}F\coloneqq\mathfrak{p}_{\Lambda}\circ\left(\Sigma_{\left|V\right.}\right)^{-1},\label{cf}
\end{equation}
where $V\subseteq N\times\Lambda$ is an open subset for which $\Sigma_{\left|V\right.}$
is a diffeomorphism onto its image.\bigskip{}

Reciprocally (see also \cite{gp}), from a submersion $F:M\rightarrow\Lambda$
satisfying \eqref{fi} and \eqref{tr}, we can construct, around every
point of $M$, a neighborhood $U$ and a local complete solution $\Sigma$
of the $\Pi$-HJE. The involved subset $U$ is one for which $\left(\Pi,F\right)_{\left|U\right.}$
is a diffeomorphism onto its image, and $\Sigma$ is given by 
\begin{equation}
\Sigma=\left(\Pi_{\left|U\right.},F_{\left|U\right.}\right)^{-1}:\Pi\left(U\right)\times F\left(U\right)\rightarrow U.\label{fc}
\end{equation}

\bigskip{}

Summarizing, a complete solution gives rise to local first integrals
\textit{via} \eqref{cf}, and first integrals give rise to a local
complete solution \textit{via} \eqref{fc}.

\subsection{Integrability by quadratures on symplectic manifolds}

\label{intq}

Let $\left(M,\omega\right)$ be a symplectic manifold. Given a distribution
$\mathcal{V}\subseteq TM$ (resp. $m\in M$ and a linear subspace
$\mathcal{V}\subseteq T_{m}M$), by $\mathcal{V}^{\bot}$ we shall
denote, as usual, the \textbf{symplectic orthogonal} of $\mathcal{V}$
w.r.t. $\omega$. The following result is a slightly extension to
general dynamical systems of a result given in \cite{gp} (only valid
for Hamiltonian systems).

\begin{thm}

\label{0} Let $F:M\rightarrow\Lambda$ be a surjective submersion
and $X\in\mathfrak{X}\left(M\right)$ a vector field such that: 
\begin{description}
\item [{I.}] $\mathsf{Im}X\subseteq\mathsf{Ker}F_{*}$ (first integrals); 
\item [{II.}] $\mathsf{Ker}F_{*}\subseteq\left(\mathsf{Ker}F_{*}\right)^{\bot}$\textbf{,
}i.e. $F$ is \textbf{isotropic}; 
\item [{III.}] and $L_{X}\beta=0$, with $\beta\coloneqq\mathfrak{i}_{X}\omega$;\footnote{If $\beta=dH$ for some function $H\in C^{\infty}\left(M\right)$,
i.e. $X=X_{H}=\omega^{\sharp}\circ dH$, this point is always satisfied.} 
\end{description}
then the trajectories of $X$ can be integrated up to quadratures.

\end{thm}

\noindent \textit{Proof.} We shall proceed in four steps.

\textbf{a.} Given a point $m\in M$, consider an open neighborhood
$U$ of $m$ and a surjective submersion $\Pi:U\rightarrow\Pi\left(U\right)$
transverse to $F_{\left|U\right.}:U\rightarrow F\left(U\right)$.
(As it is well-known, such $\Pi$ can be constructed just by fixing
a coordinate chart and solving linear equations). Using the point
(I) above and the results of the last section, it is clear (shrinking
$U$ if necessary) that $\Sigma\coloneqq\left(\Pi,F_{\left|U\right.}\right)^{-1}$
(see \eqref{fc}), which can be constructed by using the Inverse Function
Theorem, is a local complete solution of the $\Pi$-HJE for $X$.
According to Theorem 3.12 of \cite{gp} (replacing there $dH$ by
$\beta$), this implies that (recall \eqref{Srelrr}) 
\begin{equation}
\mathfrak{i}_{X^{\Sigma}}\Sigma^{*}\omega=\Sigma^{*}\beta,\label{a}
\end{equation}
omitting the restrictions to $U$ of $\omega$ and $\beta$.

\textbf{b.} Using \eqref{Srelr} and the fact that $\Sigma$ is a
diffeomorphism, the point (III) is equivalent to 
\begin{equation}
L_{X^{\Sigma}}\Sigma^{*}\beta=0.\label{lsb}
\end{equation}
For each $n\in\Pi\left(U\right)$, let us define $\beta_{n}\in\Omega^{1}\left(F\left(U\right)\right)$
such that 
\begin{equation}
\left\langle \beta_{n}\left(\lambda\right),z\right\rangle =\left\langle \Sigma^{*}\beta\left(n,\lambda\right),\left(0,z\right)\right\rangle ,\;\;\;\forall\lambda\in F\left(U\right),\;\;\;z\in T_{\lambda}\Lambda.\label{dbn}
\end{equation}
Then, along an integral curve $\left(\gamma\left(t\right),\lambda\right)$
of $X^{\Sigma}$, it can be shown from \eqref{lsb} that 
\begin{equation}
\frac{d}{dt}\left\langle \beta_{\gamma\left(t\right)}\left(\lambda\right),z\right\rangle =0\label{dtb}
\end{equation}
and, consequently, 
\[
\left\langle \beta_{\gamma\left(t\right)}\left(\lambda\right),z\right\rangle =\left\langle \beta_{\gamma\left(0\right)}\left(\lambda\right),z\right\rangle ,\;\;\;\forall\lambda\in F\left(U\right),\;\;\;z\in T_{\lambda}\Lambda.
\]
Let us prove it. From (I) and (II) we have that $\mathfrak{i}_{Y}\beta=\mathfrak{i}_{Y}\mathfrak{i}_{X}\omega=0$
for all $Y$ such that $\mathsf{Im}Y\subseteq\mathsf{Ker}F_{*}$.
This implies that 
\[
\mathfrak{i}_{\hat{Y}}\Sigma^{*}\beta=0
\]
for all $\hat{Y}\in\mathfrak{X}\left(\Pi\left(U\right)\times F\left(U\right)\right)$
of the form $\left(y,0\right)$, i.e. $\mathsf{Im}\hat{Y}\subseteq T\Pi\left(U\right)\times\left\{ 0\right\} $
(compare with \eqref{xsin}). On the other hand, given a vector field
$z\in\mathfrak{X}\left(F\left(U\right)\right)$, for $Z=\left(0,z\right)$
it is easy to see that $\left[X^{\Sigma},Z\right]$ is of the form
$\left(y,0\right)$. Then, from that and \eqref{lsb}, 
\[
L_{X^{\Sigma}}\circ\mathsf{\mathfrak{i}}_{Z}\left(\Sigma^{*}\beta\right)=\left(\mathsf{\mathfrak{i}}_{Z}\circ L_{X^{\Sigma}}+\mathfrak{i}_{\left[X^{\Sigma},Z\right]}\right)\left(\Sigma^{*}\beta\right)=\mathsf{\mathfrak{i}}_{Z}\circ L_{X^{\Sigma}}\left(\Sigma^{*}\beta\right)=0.
\]
Hence, \eqref{dtb} follows by combining \eqref{dbn} and the last
equation.

\textbf{c.} Using that $\omega$ is closed, we can assume (without
loss of generality) that 
\[
\omega=-d\theta,\mbox{ with }\theta\in\Omega^{1}(U).
\]
Now, we can proceed as in \cite{gp} (see Section 3.3.1 in \cite{gp}).
In fact, 
\[
0=\Sigma_{\lambda}^{*}\omega=-\Sigma_{\lambda}^{*}(d\theta)=-d(\Sigma_{\lambda}^{*}\theta)
\]
and, thus, a real $C^{\infty}$-function $W_{\lambda}:\Pi(U)\to{\Bbb R}$
can be constructed up to quadratures (shrinking $U$ if needed) such
that 
\[
\Sigma_{\lambda}^{*}\theta=dW_{\lambda}.
\]
In turn, the family of functions $W_{\lambda}$'s gives rise to a
real $C^{\infty}$-function $W:\Pi(U)\times F(U)\to{\Bbb R}$ satisfying
\[
\left\langle (dW-\Sigma^{*}\theta)(n,\lambda),(y,0)\right\rangle =0,\;\;\mbox{ for }(n,\lambda)\in\Pi(U)\times F(U)\mbox{ and }y\in T_{n}\Pi(U).
\]
In particular, since $X^{\Sigma}\in\mathsf{Ker}F_{*}$, it follows
that 
\[
i_{X^{\Sigma}}\Sigma^{*}\theta=i_{X^{\Sigma}}dW.
\]
Therefore, we deduce that 
\begin{eqnarray}
i_{X^{\Sigma}}\Sigma^{*}\omega=-i_{X^{\Sigma}}d(\Sigma^{*}\theta)=-{L}_{{X^{\Sigma}}}\Sigma^{*}\theta+di_{X^{\Sigma}}(\Sigma^{*}\theta)={L}_{X^{\Sigma}}(dW-\Sigma^{*}\theta).\label{b}
\end{eqnarray}
Then, combining \eqref{a} and \eqref{b}, 
\begin{equation}
L_{X^{\Sigma}}\left(dW-\Sigma^{*}\theta\right)=\Sigma^{*}\beta.\label{lxd}
\end{equation}
As a consequence, in terms of the functions $\varphi_{\lambda}:\Pi\left(U\right)\rightarrow T_{\lambda}^{*}F\left(U\right)$,
given by 
\begin{equation}
\left\langle \varphi_{\lambda}\left(n\right),z\right\rangle =\left\langle \left(dW-\Sigma^{*}\theta\right)\left(n,\lambda\right),\left(0,z\right)\right\rangle ,\;\;\;\forall z\in T_{\lambda}F\left(U\right),\label{fit}
\end{equation}
the Equation \eqref{lxd} along an integral curve $\left(\gamma\left(t\right),\lambda\right)$
of $X^{\Sigma}$ translates to (using similar calculations as in the
previous step) 
\[
\frac{d}{dt}\varphi_{\lambda}\left(\gamma\left(t\right)\right)=\beta_{\gamma\left(0\right)}\left(\lambda\right),
\]
or equivalently, 
\begin{equation}
\varphi_{\lambda}\left(\gamma\left(t\right)\right)=\varphi_{\lambda}\left(\gamma\left(0\right)\right)+t\,\beta_{\gamma\left(0\right)}\left(\lambda\right).\label{eg}
\end{equation}

\textbf{d.} Finally, since each $\varphi_{\lambda}$ is an immersion
(see Proposition 3.16, \cite{gp}), from the above equation we can
construct the curves $\gamma$ (by using the Implicit Function Theorem),
from which all the integral curves of $X_{\left|U\right.}$ can be
obtained. In fact, the latter are given by the formula $\Gamma\left(t\right)=\Sigma\left(\gamma\left(t\right),\lambda\right)$,
as explained at the end of Section \ref{comp} (see \eqref{G}). Since
all that can be done around every $m\in M$, then all the integral
curves of $X$ can be constructed in the same way. $\;\;\;\blacksquare$

\bigskip{}

Given a surjective submersion $G:M\rightarrow\Upsilon$ and a $1$-form
$\phi\in\Omega^{1}\left(\Upsilon\right)$, the vector field 
\begin{equation}
X=\omega^{\sharp}\circ G^{*}\phi\label{xwf}
\end{equation}
satisfies the point (I) above, for another submersion $F:M\rightarrow\Lambda$,
if and only if 
\begin{equation}
F_{*,m}\circ\omega^{\sharp}\circ G_{m}^{*}\left(\phi\left(G\left(m\right)\right)\right)=0,\;\;\;\forall m\in M.\label{casi}
\end{equation}
If in addition 
\begin{equation}
\mathsf{Ker}F_{*}\subseteq\mathsf{Ker}G_{*},\label{FG}
\end{equation}
then $\mathfrak{i}_{X}\circ F^{*}=\mathfrak{i}_{X}\circ G^{*}=0$
and consequently 
\begin{equation}
L_{X}\beta=L_{X}G^{*}\phi=\left(\mathfrak{i}_{X}\circ G^{*}d\phi+d\left(\mathfrak{i}_{X}\circ G^{*}\phi\right)\right)=0.\label{p3}
\end{equation}
So, given a symplectic manifold $\left(M,\omega\right)$, examples
of dynamical systems satisfying the points (I)-(III) of Theorem \ref{0}
are given by submersions $F$ and $G$ satisfying \eqref{FG}, being
$F$ isotropic, $1$-forms $\phi$ satisfying \eqref{casi} and vector
fields given by \eqref{xwf}. These examples can be seen as a generalization
of the \textit{non-commutative integrable systems}, as we show below,
and they will appear in the last section of the paper.

\subsection{Non-commutative integrability and Casimir $1$-forms}

\label{ncis}

A \textbf{Mishchenko-Fomenko} or\textbf{ non-commutative integrable
(NCI) system} \cite{mf} (see also \cite{j} and references therein)
can be defined as a triple given by a symplectic manifold $\left(M,\omega\right)$,
a Hamiltonian vector field $X_{H}=\omega^{\sharp}\circ dH$ and a
surjective submersion $F:M\rightarrow\Lambda$ such that: 
\begin{enumerate}
\item $\mathsf{Im}X_{H}\subseteq\mathsf{Ker}F_{*}$; 
\item $\mathsf{Ker}F_{*}\subseteq\left(\mathsf{Ker}F_{*}\right)^{\bot}$; 
\item $\left(\mathsf{Ker}F_{*}\right)^{\bot}$ is integrable. 
\end{enumerate}
When $\mathsf{Ker}F_{*}=\left(\mathsf{Ker}F_{*}\right)^{\bot}$, i.e.
$F$ is \textbf{Lagrangian}, the third point is automatic. In such
a case, we have a so-called \textbf{Liouville-Arnold }or\textbf{ commutative
integrable (CI) system} \cite{ar,lio}. It is well-known that all
these systems are integrable by quadratures. The traditional way of
proving that relays on the \textit{Lie theorem on integrability by
quadratures} \cite{akn,lie} (see also \cite{g}).

Usually, in the definition of NCI and CI systems, one more requirement
is included: $F$ has compact and connected leaves. In such a case,
beside integrability by quadratures, \emph{action-angle} coordinates
can be found for such systems (see \cite{rui} and \cite{eva}). We
do not consider this case here. \begin{rem} \label{ncij} An alternative
definition can be given in terms of subsets of functions $\mathcal{F}\subseteq C^{\infty}\left(M\right)$.
The conditions $3$ and $2$ above are replaced by asking $\mathcal{F}$
to be a Poisson sub-algebra and \textit{complete} (see \cite{j}),
respectively, and $1$ is replaced by asking that the elements of
$\mathcal{F}$ Poisson commute with $H$. \end{rem} To analyze the
relationship between NCI systems and the systems given at the end
of the last section, let us consider an arbitrary surjective submersion
$F:M\rightarrow\Lambda$. On the one hand, it can be shown that a
Hamiltonian vector field $X_{H}$ belongs to $\left(\mathsf{Ker}F_{*}\right)^{\bot}$
if and only if there exists a function $h:\Lambda\rightarrow\mathbb{R}$
such that $h\circ F=H$. So, if we ask that $\mathsf{Im}X_{H}\subseteq\mathsf{Ker}F_{*}$
and that $F$ is isotropic, then 
\begin{equation}
X_{H}=\omega^{\sharp}\circ F^{*}dh.\label{xhw}
\end{equation}
On the other hand, it is well-known (see, for instance, \cite{Vaisman}
Prop. 7.18) that $\left(\mathsf{Ker}F_{*}\right)^{\bot}$ is integrable
if and only if $\Lambda$ is a Poisson manifold with bi-vector $\Xi$,
given by the formula 
\[
\Xi^{\sharp}\left(\alpha\right)=F_{*,m}\circ\omega^{\sharp}\circ F_{m}^{*}\left(\alpha\right),\;\;\;\alpha\in T_{F\left(m\right)}^{*}\Lambda,
\]
and $F$ is a Poisson morphism. In such a case, the condition $\mathsf{Im}X_{H}\subseteq\mathsf{Ker}F_{*}$,
for $X_{H}$ given by \eqref{xhw}, is equivalent to (compare to \eqref{casi})
\[
\Xi^{\sharp}\left(dh\left(F\left(m\right)\right)\right)=F_{*,m}\circ\omega^{\sharp}\circ F_{m}^{*}\left(dh\left(F\left(m\right)\right)\right)=0,
\]
which says precisely that $dh$ is a \textbf{Casimir} $1$-form for
$\Xi$ (and $h$ is a Casimir function). In the case when $F$ is
Lagrangian, then $\Xi=0$, and consequently every $1$-form on $\Lambda$
is a Casimir. Thus, the NCI systems are a subclass of the systems
given at the end of the last section, where $\Lambda$ is a Poisson
manifold, $G=F:M\rightarrow\Lambda$ is a Poisson morphism and $\phi=dh$
is an exact Casimir $1$-form with respect to the Poisson structure
on $\Lambda$.

\section{Complete solutions and symmetries}

Given a general action $\rho:G\times M\rightarrow M$ (not necessarily
free or proper) of a Lie group $G$ on a manifold $M$ and a $G$-invariant
vector field $X$ on $M$ (with respect to $\rho$), we shall construct
in this section, based on $\rho$ and the canonical projection $\pi:M\rightarrow M/G$,
local fibrations $\Pi$ of $M$ and related complete solutions of
the $\Pi$-HJE for $X$. Let us begin with the local fibrations $\Pi$
based on $\pi$.

\subsection{The \textit{vertical} submersions}

\subsubsection{General actions and regular points}

Let $\rho:G\times M\rightarrow M$ be an action of a Lie group $G$
on $M$. Let us introduce some basic notation and recall some well-known
facts.

As usual, given $g\in G$ and $m\in M$, by $\rho_{g}$ and $\rho_{m}$
we shall denote the maps $\rho_{g}:M\rightarrow M$ and $\rho_{m}:G\rightarrow M$
such that $\rho_{g}\left(m\right)=\rho_{m}\left(g\right)=\rho\left(g,m\right)$.
Also, we shall denote by $\mathfrak{g}$ the Lie algebra of $G$ and
by $G_{m}$ the isotropy subgroup of $m$.

For latter convenience, let us note that 
\begin{equation}
\mathsf{Ker}\left(\rho_{m}\right)_{*,e}=\mathfrak{g}_{m},\label{kroge}
\end{equation}
where $e\in G$ is the identity element and $\mathfrak{g}_{m}$ is
the Lie algebra of $G_{m}$. And recall that the fundamental vector
field associated to $\eta\in\mathfrak{g}$ is given by 
\begin{equation}
\eta_{M}\left(m\right)\coloneqq\left(\rho_{m}\right)_{*,e}\left(\eta\right).\label{fvf}
\end{equation}

Let $\pi:M\rightarrow\left.M\right/G$ be the canonical projection
and consider on $\left.M\right/G$ the quotient topology. Recall that,
since each $\rho_{g}:M\rightarrow M$ is a homeomorphism for all $g\in G$,
then $\pi$ is open (besides continuous). Recall also the identities
\begin{equation}
\pi\circ\rho_{m}\left(g\right)=\pi\circ\rho_{g}\left(m\right)=\pi\left(m\right),\;\;\;\forall m\in M,\;g\in G.\label{piro}
\end{equation}

When $\rho$ is free (i.e. if $G_{m}=\left\{ e\right\} $ for all
$m\in M$) and proper, then, as it is well-known (see \cite{am}),
$\left.M\right/G$ has a unique manifold structure such that $\pi:M\rightarrow\left.M\right/G$
is a surjective submersion. For more general actions we shall show
a similar result, but at a local level around a \textit{regular point}.

\begin{defn} We shall say that $m_{0}\in M$ is $\rho$-\textbf{regular}
if there exists an open neighborhood $U$ of $m_{0}$ such that 
\begin{equation}
\dim G_{m}=\dim G_{m_{0}},\;\;\mbox{ for every }m\in U.\label{regcon}
\end{equation}
We shall call such neighborhood $U$ \textbf{admissible} if in addition
$U$ is connected. The (open) subset of all the $\rho$-regular points
will be denoted $\mathcal{R}_{\rho}$. \end{defn}

\begin{rem} \label{rem1} Note that if $m_{0}$ is a $\rho$-{regular}
point, then the assigning 
\[
m\in U\longmapsto{\frak{g}}_{m}\subseteq{\frak{g}}
\]
defines a vector subbundle of the trivial vector bundle $pr_{1}:U\times{\frak{g}}\to U$
for each admissible neighborhood $U$. \end{rem}

Given $m_{0}\in\mathcal{R}_{\rho}$, there exists an admissible neighborhood
$U$ of $m_{0}$ such that $\rho_{g}\left(U\right)\subseteq U$ for
all $g\in G$, i.e. $U$ is a $G$-invariant subset. To show it, note
that given $m,m'\in M$ such that $m'=\rho\left(g,m\right)$ for some
$g$, we have that $g\cdot G_{m}\cdot g^{-1}=G_{m'}$, and consequently
\[
\dim G_{m}=\dim G_{\rho_{g}\left(m\right)},\;\;\;\forall g\in G,\;m\in M.
\]
Then, given any admissible neighborhood $V$ of $m_{0}$, it is clear
that 
\begin{equation}
V_{G}\coloneqq\bigcup_{g\in G}\rho_{g}\left(V\right)\label{ug}
\end{equation}
includes $m_{0}$, is open, $G$-invariant and admissible. As a consequence,
the set $\mathcal{R}_{\rho}$ is $G$-invariant.\bigskip{}

If the action $\rho$ is free, then every element of $M$ is $\rho$-regular
and $M$ (if connected) is an admissible neighborhood. For $G=SO\left(3\right)$
acting on $M=\mathbb{R}^{3}$ with the natural action $\rho_{nat}$,
we have that $\dim G_{m}=1$ for $m\in\mathbb{R}^{3}-\left\{ 0\right\} $
and $\dim G_{0}=\dim G=3$. Thus, all the points of $\mathbb{R}^{3}$
except $0$ are $\rho_{nat}$-regular. In general, we have the following
result.

\begin{prop} $\mathcal{R}_{\rho}$ is a $G$-invariant open dense
subset of $M$. \end{prop}

\noindent \textit{Proof. } We already saw that $\mathcal{R}_{\rho}$
is $G$-invariant. We shall prove that: 
\begin{enumerate}
\item if $k$ is the minimum dimension of the isotropy subgroups and $\dim G_{m_{0}}=k$,
then $m_{0}$ is a $\rho$-regular point; 
\item the complement of $\mathcal{R}_{\rho}$ has empty interior. 
\end{enumerate}
\noindent For the first point, define 
\[
k\coloneqq\min\left\{ \dim G_{m}\;:\;m\in M\right\} 
\]
and take $m_{0}$ such that $\dim G_{m_{0}}=k$. Consider the Lie
algebra $\mathfrak{g}_{m_{0}}$ of $G_{m_{0}}$ and a linear complement
$\mathfrak{g}_{m_{0}}^{c}$ of it. For any element $v\in\mathfrak{g}_{m_{0}}^{c}-\left\{ 0\right\} $
we have that 
\[
\tilde{\rho}\left(v,m_{0}\right)\neq0,
\]
where $\tilde{\rho}$ is the action of $\mathfrak{g}$ on $M$ induced
by $\rho$. Then, by continuity, there exists a neighborhood $U$
of $m_{0}$ such that 
\[
\tilde{\rho}\left(v,m\right)\neq0,\;\;\;\forall m\in U.
\]
This means that $\dim\mathfrak{g}_{m}\leq\dim\mathfrak{g}_{m_{0}}=k$
for all $m\in U$. But $k$ is the minimum dimension, hence $\dim\mathfrak{g}_{m}=k$
for all $m\in U$. This says precisely that $m_{0}$ is a $\rho$-regular
point.\bigskip{}

For the second point, suppose that the complement $\mathcal{R}_{\rho}^{c}$
has interior, i.e., for some $m_{1}\in\mathcal{R}_{\rho}^{c}$ there
exists an open subset $U_{1}$ such that $m_{1}\in U_{1}\subseteq\mathcal{R}_{\rho}^{c}$.
Consider the Lie algebra $\mathfrak{g}_{m_{1}}$ of $G_{m_{1}}$ and
a linear complement $\mathfrak{g}_{m_{1}}^{c}$ of it. Proceeding
as above, we can ensure that there exists a neighborhood $U_{2}\subseteq U_{1}$
of $m_{1}$ such that $\dim\mathfrak{g}_{m}\leq\dim\mathfrak{g}_{m_{1}}$
for all $m\in U_{2}$. Since $U_{2}\subseteq\mathcal{R}_{\rho}^{c}$,
there must exists $m_{2}\in U_{2}$ such that $\dim\mathfrak{g}_{m_{2}}<\dim\mathfrak{g}_{m_{1}}$.
Otherwise, $m_{1}$ would be a $\rho$-regular point (with admissible
neighborhood $U_{2}$). Repeating this reasoning for $m_{2}$, we
can ensure the existence of a neighborhood $U_{3}\subseteq U_{1}$
of $m_{2}$ for which $\dim\mathfrak{g}_{m}\leq\dim\mathfrak{g}_{m_{2}}$
for all $m\in U_{3}$, and consequently the existence of a point $m_{3}\in U_{3}$
such that $\dim\mathfrak{g}_{m_{3}}<\dim\mathfrak{g}_{m_{2}}$. In
this way, since the dimension of $\mathfrak{g}$ is finite, in some
step of this procedure we shall find $m_{0}\in U_{1}\subseteq\mathcal{R}_{\rho}^{c}$
such that $\dim\mathfrak{g}_{m_{0}}=k$. Since such $m_{0}$ must
belong to $\mathcal{R}_{\rho}$, we have arrived at a contradiction.
$\;\;\;\blacksquare$

\subsubsection{The submersions $\pi_{\left|U\right.}$}

Now, let us construct smooth local versions of the canonical projection
$\pi$.

\begin{prop} Given $m_{0}\in\mathcal{R}_{\rho}$, there exists a
neighborhood $U$ for $m_{0}$ such that the open subset $\pi\left(U\right)$
has a manifold structure and the restriction $\pi_{\left|U\right.}:U\rightarrow\pi\left(U\right)$
is a submersion satisfying 
\begin{equation}
\mathsf{Ker}\left(\pi_{\left|U\right.}\right)_{*,m}=\mathsf{Im}\left(\rho_{m}\right)_{*,e},\;\;\;\forall m\in U,\label{kerim}
\end{equation}
and, consequently, 
\begin{equation}
\dim\left(\mathsf{Ker}\left(\pi_{\left|U\right.}\right)_{*}\right)=\dim G-\dim G_{m_{0}}.\label{kerpi}
\end{equation}
Moreover, $U$ can be taken $G$-invariant. \end{prop}

\noindent \textit{Proof.} Let $U_{1}$ be an admissible neighborhood
of $m_{0}$ and consider the distribution given by 
\[
\mathfrak{F}_{1}(m)\coloneqq\mathsf{Im}(\rho_{m})_{*,e}.
\]

\noindent Since $\mathfrak{F}_{1}$ is clearly generated by the fundamental
vector fields $\eta_{M}$ (see Eq. \eqref{fvf}), with $\eta\in\mathfrak{g}$,
then $\mathfrak{F}_{1}$ is involutive (see for instance \cite{am}).
And from the same reason, 
\[
\dim\left(\mathsf{Im}(\rho_{m})_{*,e}\right)=\dim\mathfrak{g}-\dim\left(\mathsf{Ker}(\rho_{m})_{*,e}\right)=\dim\mathfrak{g}-\dim\mathfrak{g}_{m}=\dim G-\dim G_{m},
\]
which is constant and equal to $r_{1}\coloneqq\dim G-\dim G_{m_{0}}$
for all $m\in U_{1}$ (because of \eqref{regcon}). Then, defining
$r\coloneqq\dim M$ and using the Frobenius Theorem, we can find a
local chart 
\[
(U_{2},\varphi\equiv(x^{1},\dots,x^{r_{1}},x^{r_{1}+1},\dots,x^{r}))
\]
in $U_{1}$ such that $m_{0}\in U_{2}$, 
\[
\varphi(U_{2})=V_{2}\times V_{2}'\subseteq{\Bbb R}^{r_{1}}\times{\Bbb R}^{r-r_{1}},
\]
with $V_{2}$ and $V_{2}'$ open subsets in ${\Bbb R}^{r_{1}}$ and
${\Bbb R}^{r-r_{1}}$, respectively, and 
\[
\mathfrak{F}_{1}(m)=\left\langle \left.\frac{\partial}{\partial x^{1}}\right|_{m},\dots,\left.\frac{\partial}{\partial x^{r_{1}}}\right|_{m}\right\rangle \mbox{ for all }m\in U_{2}.
\]
Now, we can consider the $G$-invariant open subset $U$ of $M$ given
by 
\[
U=\bigcup_{g\in G}\rho_{g}(U_{2}).
\]
It is clear that $m_{0}\in U_{2}\subseteq U$ and, moreover, $U/G\cong V_{2}'$
and the canonical projection $\pi_{|U}:U\to U/G\cong V_{2}'$ is a
surjective submersion. $\;\;\;\blacksquare$

\bigskip{}

From now on, by \textbf{\textit{admissible}} we shall mean any admissible
neighborhood $U$ of $m_{0}$ for which the last proposition holds.

\bigskip{}

The following result will be useful later.

\begin{prop}\label{prop5'} Given $m_{0}\in\mathcal{R}_{\rho}$ and
an admissible neighborhood $U$ of $m_{0}$, the subset 
\begin{equation}
R_{U}\coloneqq\left\{ \left(m,\rho_{g}\left(m\right)\right)\;:\;g\in G,\;m\in U\right\} \subseteq U\times M\label{ru}
\end{equation}
is a closed submanifold of $U\times M$, of dimension $\dim M+\dim G-\dim G_{m_{0}}$,
and the surjective map 
\begin{equation}
\Phi_{U}:\left(g,m\right)\in G\times U\longmapsto\left(m,\rho_{g}\left(m\right)\right)\in R_{U}\label{Phi-U}
\end{equation}
is smooth and also a submersion. \end{prop}

\textit{Proof. } Consider the admissible $G$-invariant open subset
$\tilde{U}=\bigcup_{g\in G}\rho_{g}\left(U\right)$ (see Eq. \eqref{ug}),
the related subset $R_{\tilde{U}}$ and the related surjective map
$\Phi_{\tilde{U}}:G\times\tilde{U}\to R_{\tilde{U}}$ (given as in
\eqref{ru} and \eqref{Phi-U}). If we prove the proposition for $\tilde{U}$,
since $R_{U}=R_{\tilde{U}}\cap\left(U\times M\right)$ and $\Phi_{U}=\left.\Phi_{\tilde{U}}\right|_{G\times U}$,
then we would proved it for $U$.

\noindent Using that the space of orbits $\pi(\tilde{U})=\tilde{U}/G$
is a quotient smooth manifold and a classical result (see, for instance,
\cite{Lee}), we deduce that $R_{\tilde{U}}\subseteq\tilde{U}\times\tilde{U}$
is a closed submanifold of $\tilde{U}\times\tilde{U}$ (and also of
$\tilde{U}\times M$). As a consequence, since 
\[
\left(g,m\right)\in G\times\tilde{U}\longmapsto\left(m,\rho_{g}\left(m\right)\right)\in\tilde{U}\times\tilde{U}
\]
is smooth, the same is true for the surjection $\Phi_{\tilde{U}}:G\times\tilde{U}\rightarrow R_{\tilde{U}}$.
To find the dimension of $R_{\tilde{U}}$ and show that $\Phi_{\tilde{U}}$
is a submersion, it is enough to calculate the rank of $\Phi_{\tilde{U}}$
and show that is constant (since $\Phi_{\tilde{U}}$ is surjective).
Let us do that.

From (\ref{Phi-U}) and the identity $\rho_{m}\circ L_{g}=\rho_{g}\circ\rho_{m}$,
it follows that 
\begin{equation}
(\Phi_{\tilde{U}})_{*,(g,m)}(u,v)=(v,(\rho_{g})_{*,m}(v+(\rho_{m})_{*,e}((L_{g^{-1}})_{*,g}(u)))),\quad\forall u\in T_{g}G,\;v\in T_{m}\tilde{U},\label{Phi-g-m-0}
\end{equation}
and in particular, for $g=e$, 
\begin{equation}
(\Phi_{\tilde{U}})_{*,(e,m)}(\eta,v)=(v,v+(\rho_{m})_{*,e}(\eta)),\quad\forall\eta\in\mathfrak{g},\;v\in T_{m}\tilde{U}.\label{Phi-e-m}
\end{equation}
Then, from \eqref{Phi-g-m-0} and \eqref{Phi-e-m} we have that 
\[
(\Phi_{\tilde{U}})_{*,(g,m)}=(id_{T_{m}M}\times(\rho_{g})_{*,m})\circ(\Phi_{\tilde{U}})_{*,(e.m)}\circ((L_{g^{-1}})_{*,g}\times id_{T_{m}M}).
\]
Consequently, for all $\left(g,m\right)\in G\times\tilde{U}$, 
\[
(\mbox{rank }\Phi_{\tilde{U}})(g,m)=(\mbox{rank }\Phi_{\tilde{U}})(e,m)=\dim M+(\dim G-\dim G_{m_{0}}),
\]
which ends our proof. $\;\;\;\blacksquare$

\subsubsection{Symplectic actions and momentum maps}

Suppose that $M$ is a symplectic manifold, with symplectic form $\omega$,
and $\rho$ is a \textbf{symplectic action}, i.e. 
\[
\left(\rho_{g}\right)^{*}\circ\omega=\omega,\;\;\;\forall g\in G,
\]
or equivalently 
\begin{equation}
\left(\rho_{g}\right)_{*,m}\circ\omega_{m}^{\sharp}=\omega_{\rho_{g}\left(m\right)}^{\sharp}\circ\left(\rho_{g^{-1}}\right)_{\rho_{g}\left(m\right)}^{*},\;\;\;\forall m\in M,\;g\in G.\label{sya}
\end{equation}

\begin{prop}\label{pm} Under the above conditions, for each admissible
neighborhood $U$ we have that: 
\begin{enumerate}
\item The manifold $\pi\left(U\right)$ has a Poisson structure $\Xi_{U}$,
characterized by the condition 
\begin{equation}
\Xi_{U}\left(\alpha,\beta\right)\circ\pi=\omega\left(\omega^{\sharp}\left(\pi^{*}\alpha\right),\omega^{\sharp}\left(\pi^{*}\beta\right)\right),\;\;\;\forall\alpha,\beta\in\Omega^{1}(\pi(U)),\label{T}
\end{equation}
with respect to which $\pi_{\left|U\right.}$ is a Poisson morphism. 
\item Let $X$ be a $G$-invariant vector field, i.e. 
\begin{equation}
\left(\rho_{g}\right)_{*}\circ X=X\circ\rho_{g},\;\;\;\forall g\in G.\label{ginv}
\end{equation}
Then there exists a unique vector field $Y\in\mathfrak{X}\left(\pi\left(U\right)\right)$
such that 
\[
\left(\pi_{\left|U\right.}\right)_{*}\circ X_{\left|U\right.}=Y\circ\pi_{\left|U\right.}.
\]
\end{enumerate}
\end{prop}

\noindent \textit{Proof. } $\left(1\right)$ This result is proved
in \cite{mr} under the hypothesis that $U$ is $G$-invariant and
that the $G$-action on $U$ is free and proper. But, in that proof
the key point is that the space of orbits $\pi\left(U\right)$ is
a quotient manifold, as in our case.

\medskip{}

\noindent $(2)$ It is also a well-known result (see, for instance
\cite{Quam}) that if $U$ is a principal $G$-bundle over $U/G$,
then every $G$-invariant vector field over $U$ is projectable over
$U/G$. But, as in $(1)$, the key point in order to prove this fact
is that $U/G$ is a quotient manifold. So, proceeding in a similar
way as in \cite{Quam}, we deduce $(2)$. $\;\;\;\blacksquare$

\bigskip{}

Suppose that $\rho$ has a (global) momentum map, i.e. a function
$K:M\rightarrow\mathfrak{g}^{*}$ such that 
\begin{equation}
\left\langle \omega^{\flat}\left(\left(\rho_{m}\right)_{*,e}\left(\eta\right)\right),v\right\rangle =\left\langle K_{*,m}\left(v\right),\eta\right\rangle ,\;\;\;\forall m\in M,\;v\in T_{m}M,\;\eta\in\mathfrak{g}.\label{defk}
\end{equation}

\begin{prop} For each admissible neighborhood $U$, 
\begin{equation}
\mathsf{Ker}\left(K_{\left|U\right.}\right)_{*}=\left(\mathsf{Ker}\left(\pi_{\left|U\right.}\right)_{*}\right)^{\bot}.\label{kkkp}
\end{equation}

\end{prop}

For a proof of this result see, for instance, \cite{am}.

\bigskip{}

Suppose in addition that $K$ can be chosen $Ad^{*}$-equivariant,
i.e. 
\begin{equation}
K\left(\rho\left(g,m\right)\right)=Ad_{g}^{*}\left(K\left(m\right)\right).\label{adinv}
\end{equation}
Here, as usual, $Ad:G\times\mathfrak{g}\rightarrow\mathfrak{g}:\left(g,\eta\right)\mapsto Ad_{g}\eta$
denotes the \textbf{adjoint action} and 
\[
Ad^{*}:G\times\mathfrak{g}^{*}\rightarrow\mathfrak{g}^{*}:\left(g,\alpha\right)\longmapsto\left(Ad_{g^{-1}}\right)^{*}\left(\alpha\right)
\]
the \textbf{co-adjoint} one.

\begin{prop} \label{5} If $K\left(\mathcal{R}_{\rho}\right)\cap\mathcal{R}_{Ad^{*}}\neq\emptyset$,
then there exists a $G$-invariant open subset $U\subseteq\mathcal{R}_{\rho}$
such that 
\[
K\left(U\right)\subseteq\mathcal{R}_{Ad^{*}}.
\]
\end{prop}

\noindent \textit{Proof. } Let $m_{1}\in\mathcal{R}_{\rho}$ be such
that $K\left(m_{1}\right)\in\mathcal{R}_{Ad^{*}}\subseteq\mathfrak{g}^{*}$,
and let $U_{1}$ be an admissible neighborhood of $m_{1}$. Given
a $G$-invariant admissible neighborhood $V\subseteq\mathcal{R}_{Ad^{*}}$
of $K\left(m_{1}\right)$ (with respect to the co-adjoint action),
define 
\[
U\coloneqq\bigcup_{g\in G}\rho_{g}\left(U_{1}\cap K^{-1}\left(V\right)\right).
\]
It is clear that $U$ is a $G$-invariant open subset and $U\subseteq\mathcal{R}_{\rho}.$
Moreover, $K\left(m\right)\in V$ for all $m\in U_{1}\cap K^{-1}\left(V\right)$.
Then $Ad_{g}^{*}K\left(m\right)\in V$ for all $g\in G$, because
of the $G$-invariance of $V$, and consequently (see \eqref{adinv})
\[
K\left(\rho_{g}\left(m\right)\right)=Ad_{g}^{*}K\left(m\right)\in V.
\]
This completes our proof. $\;\;\;\blacksquare$

\bigskip{}

The previous result will be useful at the end of the paper.

\subsection{The horizontal submersions}

In this subsection, for each admissible neighborhood $U$, we shall
construct submersions $\Theta$ transverse to the restricted canonical
projection $\pi_{\left|U\right.}$. In terms of such submersions $\Theta$,
we shall present at the end of the section the complete solutions
we are looking for.

\subsubsection{Trivializations and (local) flat connections for principal bundles}

\label{tlf}

Suppose that $\rho:G\times M\rightarrow M$ is free and proper and
consider the associated principal $G$-bundle $\pi:M\rightarrow\left.M\right/G$.
Given a local section $s:V\rightarrow U$ of $\pi$, with $U=\pi^{-1}\left(V\right)\subseteq M$
and $V\subseteq\left.M\right/G$ an open subset, we have a trivialization
\begin{equation}
\Psi=\left(\pi,\psi\right):U\rightarrow\pi\left(U\right)\times G\label{trivfi}
\end{equation}
given by $\Psi\left(\rho\left(g,s\left(\lambda\right)\right)\right)=\left(\lambda,g\right)$,
or equivalently 
\[
\Psi^{-1}\left(\lambda,g\right)=\rho\left(g,s\left(\lambda\right)\right),\;\;\;\forall\lambda\in\pi\left(U\right),\;g\in G.
\]
($\Psi$ is well-defined and invertible because $\rho$ is free).
Note that the map $\psi:U\rightarrow G$ satisfies 
\begin{equation}
\rho\left(\psi\left(m\right),s\left(\pi\left(m\right)\right)\right)=\Psi^{-1}\left(\pi\left(m\right),\psi\left(m\right)\right)=m,\;\;\;\forall m\in U.\label{rots0}
\end{equation}
Also, 
\[
\psi\left(\rho\left(g,m\right)\right)=\psi\left(\rho\left(g,\rho\left(\psi\left(m\right),s\left(\pi\left(m\right)\right)\right)\right)\right)=\psi\left(\rho\left(g\,\psi\left(m\right),s\left(\pi\left(m\right)\right)\right)\right)=g\,\psi\left(m\right),
\]
and consequently 
\begin{equation}
\psi\circ\rho_{g}=L_{g}\circ\psi\;\;\;\textrm{and}\;\;\;\psi\circ\rho_{m}=R_{\psi\left(m\right)}.\label{iden}
\end{equation}
On the other hand, it is easy to show that the map $A:TU\rightarrow\mathfrak{g}$
given by 
\begin{equation}
A\left(v\right)=\left(R_{\psi\left(m\right)}\right)_{*,e}^{-1}\psi_{*,m}\left(v\right),\mbox{ for all \ensuremath{v\in T_{m}M,}}\label{av}
\end{equation}
is a local principal connection for $\pi$. In fact, for all $m\in U$,
it follows from \eqref{iden} and \eqref{av} that 
\begin{equation}
A\left(\left(\rho_{m}\right)_{*,e}\left(\eta\right)\right)=\eta,\;\;\;\forall\eta\in\mathfrak{g},\label{c1}
\end{equation}
and 
\begin{equation}
A\left(\left(\rho_{g}\right)_{*,m}\left(v\right)\right)=Ad_{g}\left(A\left(v\right)\right),\;\;\;\forall v\in T_{m}U,\;\;\;\forall g\in G.\label{c2}
\end{equation}
In addition, since $\mathsf{Ker}A=\mathsf{Ker}\psi_{*}$, the horizontal
distribution is integrable, i.e. $A$ is a flat connection. In the
next section, we shall construct an object similar to $A$, but related
to an arbitrary action and its regular points.

\subsubsection{A flat-connection-like object for $\pi_{\left|U\right.}$}

Now, suppose that $\rho$ is a general Lie group action. For each
$\rho$-regular point $m_{0}$, we shall construct a family of submersions
transverse to $\pi_{\left|U\right.}$ (being $U$ an admissible neighborhood
of $m_{0}$). To do that, we need the next results.

\begin{lem} \label{l1}Let $G:P\rightarrow Q$ be a submersion, $p_{0}\in P$
and $\mathcal{W}\subseteq T_{p_{0}}P$ a linear complement of $\mathsf{Ker}G_{*,p_{0}}$.
Then, there exists a neighborhood $V$ of $G\left(p\right)\in Q$
and a local section $S:V\rightarrow P$ of $G$ such that 
\[
S\left(G\left(p_{0}\right)\right)=p_{0}\;\;\;\textrm{and}\;\;\;\mathsf{Im}S_{*,G\left(p_{0}\right)}=\mathcal{W}.
\]
\end{lem}

\noindent \textit{Proof. }Let $\varphi=\left(x_{1},...,x_{n}\right):U\rightarrow\varphi\left(U\right)$
be a coordinate system of $P$ around $p_{0}$. Consider the annihilator
$\mathcal{W}^{0}\subseteq T_{p_{0}}^{*}P$ of $\mathcal{W}$ and suppose
that the co-vectors 
\[
\xi_{i}=\sum_{j=1}^{n}w_{i}^{j}\,dx_{j}\left(p_{0}\right),\;\;\;i=1,...,k,
\]
give a basis for $\mathcal{W}^{0}$. Define $F:U\rightarrow\mathbb{R}^{k}$
as 
\[
F\left(\varphi^{-1}\left(x_{1},...,x_{n}\right)\right)=\left(\sum_{j=1}^{n}w_{1}^{j}\,x_{j},...,\sum_{j=1}^{n}w_{k}^{j}\,x_{j}\right).
\]
It is clear that $\mathsf{Ker}F_{*,p_{0}}=\mathcal{W}$. Then, since
$\mathsf{Ker}G_{*,p_{0}}$ and $\mathcal{W}$ are complementary, $\left(G,F\right)$
is a diffeomorphism onto its image $G\left(U\right)\times F\left(U\right)$,
shrinking $U$ if needed. As a consequence, the function $S:G\left(U\right)\rightarrow P$
such that 
\[
S\left(q\right)=\left(G,F\right)^{-1}\left(q,F\left(p_{0}\right)\right)
\]
is a smooth local section of $G$ and satisfies $S\left(G\left(p_{0}\right)\right)=p_{0}$.
Also, given $w\in\mathcal{W}$, $S_{*,G\left(p_{0}\right)}\left(G_{*,p_{0}}\left(w\right)\right)=w$.
In particular, since $G_{*,p_{0}}$ is surjective, even restricted
to $\mathcal{W}$, then $\mathsf{Im}S_{*,G\left(p_{0}\right)}=\mathcal{W}$.
So, the wanted result follows for $V=G\left(U\right)$. $\;\;\;\blacksquare$\bigskip{}

\noindent Note that the construction of the section $S$ has been
made just by using algebraic manipulations and the Inverse Function
Theorem.\bigskip{}

For the rest of the section, fix a $\rho$-regular point $m_{0}$,
an admissible neighborhood $U$ and a section $s:\pi\left(U\right)\rightarrow U$
of $\pi_{\left|U\right.}$ such that 
\begin{equation}
s\left(\pi\left(m_{0}\right)\right)=m_{0}.\label{spi}
\end{equation}

\begin{lem} \label{l2}The function 
\[
\mathcal{F}:G\times\pi\left(U\right)\rightarrow M:\left(g,\lambda\right)\mapsto\rho\left(g,s\left(\lambda\right)\right)
\]
is an open map around $\left(e,\pi(m_{0})\right)$. \end{lem}

\noindent \textit{Proof. } Note first that (according to \eqref{kerpi})
\[
\dim\left(\pi\left(U\right)\right)=\dim M-\left(\dim G-\dim G_{m_{0}}\right),
\]
and consequently 
\[
\dim\left(G\times\pi\left(U\right)\right)=\dim M+\dim G_{m_{0}}.
\]
So, it is enough to show that $\dim\left(\mathsf{Ker}\mathcal{F}_{*,\left(e,\pi(m_{0})\right)}\right)=\dim G_{m_{0}}$.
Given $X\in\mathfrak{g}$ and $Y\in T_{\pi\left(m_{0}\right)}\pi\left(U\right)$,
if 
\[
0=\mathcal{F}_{*,\left(e,\pi\left(m_{0}\right)\right)}\left(X,Y\right)=\rho_{*,\left(e,m_{0}\right)}\left(X,s_{*,\pi\left(m_{0}\right)}\left(Y\right)\right)=\left(\rho_{m_{0}}\right)_{*,e}\left(X\right)+\left(\rho_{e}\right)_{*,m_{0}}\left(s_{*,\pi\left(m_{0}\right)}\left(Y\right)\right)
\]
(where we have used \eqref{spi}), then, applying $\left(\pi_{\left|U\right.}\right)_{*,m_{0}}$
above, 
\[
0=\left(\pi_{\left|U\right.}\right)_{*,m_{0}}\circ\left(\rho_{m_{0}}\right)_{*,e}\left(X\right)+\left(\pi_{\left|U\right.}\right)_{*,m_{0}}\circ\left(\rho_{e}\right)_{*,m_{0}}\circ\left(s_{*,\pi\left(m_{0}\right)}\left(Y\right)\right).
\]
On the other hand, from \eqref{piro} we have that 
\begin{equation}
\left(\pi_{\left|U\right.}\right)_{*,m_{0}}\circ\left(\rho_{m_{0}}\right)_{*,e}=0\;\;\;\textrm{and}\;\;\;\left(\pi_{\left|U\right.}\right)_{*,m_{0}}\circ\left(\rho_{e}\right)_{*,m_{0}}=\left(\pi_{\left|U\right.}\right)_{*,m_{0}}.\label{piro2}
\end{equation}
Since in addition $\left(\pi_{\left|U\right.}\right)_{*,m_{0}}\circ s_{*,\pi\left(m_{0}\right)}$
is the identity, then $Y=0$. Hence, $\mathsf{Ker}\mathcal{F}_{*,\left(e,m_{0}\right)}$
is given by the vectors $\left(X,0\right)$ such that $\left(\rho_{m_{0}}\right)_{*,e}\left(X\right)=0$,
i.e. $X\in\mathfrak{g}_{m_{0}}$ (recall \eqref{kroge}). This ends
our proof. $\;\;\;\blacksquare$

\noindent \medskip{}

Now, the main result of the section.

\begin{thm} \label{pro-tita} Given an admissible neighborhood $U$
of $m_{0}\in\mathcal{R}_{\rho}$ and a section $s:\pi\left(U\right)\rightarrow U$
of $\pi_{\left|U\right.}$ satisfying \eqref{spi}, we can construct,
shrinking $U$ if necessary, a surjective submersion $\Theta:U\rightarrow\Theta\left(U\right)\subseteq G$
transverse to $\pi_{\left|U\right.}$ (see \eqref{tr}) such that
$\Theta\left(m_{0}\right)=e$ and 
\begin{equation}
\rho\left(\Theta\left(m\right),s\left(\pi\left(m\right)\right)\right)=m,\;\;\;\forall m\in U.\label{rots}
\end{equation}
We shall call \textbf{$s$-horizontal}, or simply \textbf{horizontal},
to such submersions $\Theta$.

\end{thm}

\noindent \textit{Proof. } First, let us make some observations about
the submersion $\Phi_{U}$ of Proposition \ref{prop5'}. 
\begin{itemize}
\item We have that $\Phi_{U}\left(e,m\right)=\left(m,m\right)$ and from
(\ref{Phi-e-m}) it follows that 
\begin{equation}
\mathsf{Ker}\left(\Phi_{U}\right){}_{*,\left(e,m\right)}=\mathfrak{g}_{m}\times\left\{ 0\right\} .\label{3.26'}
\end{equation}
\item Then, using Lemma \ref{l1}, we can construct a local section $S$
of the submersion $\Phi_{U}$ such that $S\left(m_{0},m_{0}\right)=\left(e,m_{0}\right)$
and 
\begin{equation}
\mathsf{Im}\left(S_{*,\left(m_{0},m_{0}\right)}\right)=\mathfrak{g}_{m_{0}}^{c}\times T_{m_{0}}M,\label{ims}
\end{equation}
being $\mathfrak{g}_{m_{0}}^{c}\subseteq\mathfrak{g}$ some complement
of $\mathfrak{g}_{m_{0}}$ (because \eqref{3.26'} and \eqref{ims}
are complementary). 
\end{itemize}
For simplicity, let us restrict $\Phi_{U}$ to a subset $W\times U'\subseteq G\times U$,
with $W\subseteq G$ an open neighborhood of $e$ and $U'\subseteq U\subseteq M$
an open neighborhood of $m_{0}$, such that the above mentioned section
$S$ becomes a global section $S:\Phi_{U}\left(W\times U'\right)\rightarrow W\times U'$.
Moreover, take $U'$ such that 
\begin{equation}
s\left(\pi\left(U'\right)\right)\subseteq U',\label{pos}
\end{equation}
which can be done because of \eqref{spi}. Let us write 
\[
S\left(m_{1},m_{2}\right)=\left(g_{S}\left(m_{1},m_{2}\right),m_{1}\right).
\]
Notice that, since 
\[
\left(m_{1},m_{2}\right)=\Phi_{U}\left(S\left(m_{1},m_{2}\right)\right)=\Phi_{U}\left(g_{S}\left(m_{1},m_{2}\right),m_{1}\right)=\left(m_{1},\rho\left(g_{S}\left(m_{1},m_{2}\right),m_{1}\right)\right),
\]
then 
\begin{equation}
\rho\left(g_{S}\left(m_{1},m_{2}\right),m_{1}\right)=m_{2}.\label{cle}
\end{equation}
On the other hand, since $S\left(m_{0},m_{0}\right)=\left(e,m_{0}\right)$,
then 
\begin{equation}
g_{S}\left(m_{0},m_{0}\right)=e.\label{gmme}
\end{equation}
Also, from \eqref{ims} it follows that $\mathsf{Im}\left(g_{S}\right)_{*,\left(m_{0},m_{0}\right)}\subseteq\mathfrak{g}_{m_{0}}^{c}$,
and consequently (recall \eqref{kroge}) 
\begin{equation}
\left(\rho_{m_{0}}\right)_{*,e}\circ\left(g_{S}\right)_{*,\left(m_{0},m_{0}\right)}\left(x,y\right)=0\;\;\;\Longleftrightarrow\;\;\;\left(g_{S}\right)_{*,\left(m_{0},m_{0}\right)}\left(x,y\right)=0.\label{ii}
\end{equation}

Now, consider the subset 
\[
U''\coloneqq U'\cap\left(\rho\left(W,s\left(\pi\left(U'\right)\right)\right)\right).
\]
According to Lemma \ref{l2}, $U''$ is open (shrinking $W$ and $U'$
if needed) and, since $m_{0}$ is there (see \eqref{spi}), it is
nonempty. Finally, define $\Theta:U''\rightarrow G$ as 
\[
\Theta\left(m\right)=g_{S}\left(s\left(\pi\left(m\right)\right),m\right).
\]
Let us see that $\Theta$ is well defined. If $m\in U''$, then $m\in U'$
and 
\[
m=\rho\left(g,s\left(\pi\left(m'\right)\right)\right),\;\;\;\textrm{with}\;g\in W\;\textrm{and}\;m'\in U'.
\]
Then, applying $\pi$ on both members of above equality, it follows
that $\pi\left(m\right)=\pi\left(m'\right)$, and consequently 
\[
m=\rho\left(g,s\left(\pi\left(m\right)\right)\right),\;\;\;\textrm{with}\;g\in W.
\]
In addition, since $m\in U'$, then $s\left(\pi\left(m\right)\right)\in U'$
(see \eqref{pos}). Thus, given $m\in U''$ we have that 
\[
\left(s\left(\pi\left(m\right)\right),m\right)=\left(s\left(\pi\left(m\right)\right),\rho\left(g,s\left(\pi\left(m\right)\right)\right)\right)\in\Phi_{U}\left(W\times U'\right),
\]
i.e. $\left(s\left(\pi\left(m\right)\right),m\right)$ belongs to
the domain of $S$. From \eqref{gmme}, it is clear that $\Theta\left(m_{0}\right)=e$
and, using \eqref{cle}, the identity \eqref{rots} follows. A direct
consequence of the latter is that, for all $m\in U''$, 
\begin{equation}
\left(id_{M}\right){}_{*,m}=\left(\rho_{\Theta\left(m\right)}\right)_{*,s\left(\pi\left(m\right)\right)}\circ s_{*,\pi\left(m\right)}\circ\left(\pi_{\left|U\right.}\right)_{*,m}+\left(\rho_{s\left(\pi\left(m\right)\right)}\right)_{*,\Theta\left(m\right)}\circ\Theta_{*,m},\label{idpitita}
\end{equation}
which in turn implies that 
\begin{equation}
\mathsf{Ker}\left(\pi_{\left|U\right.}\right)_{*,m}\cap\mathsf{Ker}\Theta_{*,m}=\left\{ 0\right\} ,\;\;\;\forall m\in U''.\label{nuc}
\end{equation}
Also, it implies that 
\begin{equation}
v\in\mathsf{Ker}\left(\pi_{\left|U\right.}\right)_{*,m}\;\;\;\Longleftrightarrow\;\;\;v=\left(\rho_{s\left(\pi\left(m\right)\right)}\right)_{*,\Theta\left(m\right)}\circ\Theta_{*,m}\left(w\right),\label{vker}
\end{equation}
for some $w$. Let us show it. The first implication is immediate
by applying both sides of \eqref{idpitita} to $v$, and it is fulfilled
for $w=v$. For the converse, it is enough to note that, from \eqref{piro},
\[
\left(\pi_{\left|U\right.}\right)_{*,m}\circ\left(\rho_{s\left(\pi\left(m\right)\right)}\right)_{*,\Theta\left(m\right)}=0.
\]
Something similar to \eqref{vker} can be said about $\mathsf{Ker}\Theta_{*,m}$
for $m=m_{0}$. Let us see that. Equation \eqref{idpitita} for $m=m_{0}$
reduces to 
\begin{equation}
\left(id_{M}\right){}_{*,m_{0}}=s_{*,\pi\left(m_{0}\right)}\circ\left(\pi_{\left|U\right.}\right)_{*,m_{0}}+\left(\rho_{m_{0}}\right)_{*,e}\circ\Theta_{*,m_{0}},\label{idpitita2}
\end{equation}
since $s\left(\pi\left(m_{0}\right)\right)=m_{0}$ and $\Theta\left(m_{0}\right)=e$.
Then 
\begin{equation}
v\in\mathsf{Ker}\Theta_{*,m_{0}}\;\;\;\Longleftrightarrow\;\;\;v=s_{*,\pi\left(m_{0}\right)}\circ\left(\pi_{\left|U\right.}\right)_{*,m_{0}}\left(w\right),\label{vhor0}
\end{equation}
for some $w$. The first implication follows by applying both sides
of \eqref{idpitita2} to $v$, and it is fulfilled for $w=v$. For
the converse, note first that, using \eqref{idpitita2}, 
\[
s_{*,\pi\left(m_{0}\right)}\circ\left(\pi_{\left|U\right.}\right)_{*,m_{0}}\left(w\right)=v=s_{*,\pi\left(m_{0}\right)}\circ\left(\pi_{\left|U\right.}\right)_{*,m_{0}}\left(v\right)+\left(\rho_{m_{0}}\right)_{*,e}\circ\Theta_{*,m_{0}}\left(v\right).
\]
Then, applying $\left(\pi_{\left|U\right.}\right)_{*,m_{0}}$ to the
first and the last members and using the first part of \eqref{piro2}
and the fact that $s$ is a section of $\pi_{\left|U\right.}$, we
have that 
\[
\left(\pi_{\left|U\right.}\right)_{*,m_{0}}\left(w\right)=\left(\pi_{\left|U\right.}\right)_{*,m_{0}}\left(v\right),
\]
and consequently $v=s_{*,\pi\left(m_{0}\right)}\circ\left(\pi_{\left|U\right.}\right)_{*,m_{0}}\left(v\right)$.
Finally, combining \eqref{ii} and \eqref{idpitita2}, the converse
of \eqref{vhor0} follows. So, from \eqref{vker} at $m_{0}$ and
\eqref{vhor0}, 
\begin{equation}
\mathsf{Ker}\left(\pi_{\left|U\right.}\right)_{*,m_{0}}+\mathsf{Ker}\Theta_{*,m_{0}}=T_{m_{0}}M.\label{sum}
\end{equation}
As a consequence (from \eqref{nuc} and \eqref{sum}), there exists
an admissible neighborhood $\hat{U}\subseteq U''$ of $m_{0}$ such
that 
\begin{equation}
\mathsf{Ker}\left(\pi_{\left|\hat{U}\right.}\right)_{*,m}\oplus\mathsf{Ker}\Theta_{*,m}=T_{m}M,\;\;\;\forall m\in\hat{U},\label{sud}
\end{equation}
what tell us that the rank of $\Theta$ is constant and given by $\mathsf{k}=\dim G-\dim G_{m_{0}}$
(see \eqref{kerpi}). In resume, using the Constant Rank Theorem,
we can say that, shrinking the original admissible neighborhood $U$
(if necessary), $\Theta\left(U\right)\subseteq G$ is a closed $\mathsf{k}$-dimensional
submanifold and $\Theta:U\rightarrow\Theta\left(U\right)$ is a surjective
submersion transverse to $\pi_{\left|U\right.}$, as we wanted to
show. $\;\;\;\blacksquare$

\medskip{}

\begin{rem*} It is worth mentioning that, combining \eqref{idpitita},
\eqref{vker} and \eqref{sud}, it follows for all $m\in U$ that
\begin{equation}
v\in\mathsf{Ker}\Theta_{*,m}\;\;\;\Longleftrightarrow\;\;\;v=\left(\rho_{\Theta\left(m\right)}\right)_{*,s\left(\pi\left(m\right)\right)}\circ s_{*,\pi\left(m\right)}\circ\left(\pi_{\left|U\right.}\right)_{*,m}\left(v\right).\label{vkt}
\end{equation}
\end{rem*}

Let us study some properties of $\Theta$.

\begin{prop} For any $s$-horizontal submersion $\Theta:U\rightarrow\Theta\left(U\right)$
we have that 
\begin{equation}
\Theta\left(s\left(\lambda\right)\right)=e,\;\;\;\forall\lambda\in\pi\left(U\right).\label{tsl}
\end{equation}
\end{prop}

\noindent \textit{Proof. } First, recall that $U$ is connected (ipso
facto path connected) and consequently the same is true for $\pi\left(U\right)$.
On the one hand, given $\lambda\in\Pi\left(U\right)$, it follows
from \eqref{rots} for $m=s\left(\lambda\right)$ that 
\[
\rho\left(\Theta\left(s\left(\lambda\right)\right),s\left(\lambda\right)\right)=s\left(\lambda\right),
\]
what implies that $\Theta\left(s\left(\lambda\right)\right)\in G_{s\left(\lambda\right)}$.
On the other hand, it is easy to see that, for all $g\in G$ and $m\in M$,
\[
T_{g}G_{m}=\mathsf{Ker}\left(\rho_{m}\right)_{*,g}.
\]
Then, for every vector $w\in T_{\lambda}\pi\left(U\right)$, 
\[
\Theta_{*,s\left(\lambda\right)}\left(s_{*,\lambda}\left(w\right)\right)\in\mathsf{Ker}\left(\rho_{s\left(\lambda\right)}\right)_{*,\Theta\left(s\left(\lambda\right)\right)}.
\]
As a consequence, applying \eqref{idpitita} to $v=s_{*,\lambda}\left(w\right)$
and using \eqref{vkt}, we have that 
\[
s_{*,\lambda}\left(w\right)\in\mathsf{Ker}\Theta_{*,s\left(\lambda\right)}.
\]
So, given a curve $t\in\left(-\epsilon,\epsilon\right)\mapsto\lambda\left(t\right)\in\pi\left(U\right)$
such that $\lambda\left(0\right)=\pi\left(m_{0}\right)$, we have
that 
\[
\Theta\left(s\left(\lambda\left(0\right)\right)\right)=\Theta\left(s\left(\pi\left(m_{0}\right)\right)\right)=\Theta\left(m_{0}\right)=e
\]
and 
\[
\frac{d}{dt}\Theta\left(s\left(\lambda\left(t\right)\right)\right)=0,\;\;\;\forall t\in\left(-\epsilon,\epsilon\right),
\]
from which, and the fact that $\pi\left(U\right)$ is connected, the
proposition follows. $\;\;\;\blacksquare$

\bigskip{}

Also, given $m\in U$ and $g\in G$ such that $\rho\left(g,m\right)\in U$,
it can be shown from \eqref{rots} that 
\begin{equation}
\Theta\left(\rho\left(g,m\right)\right)=g\cdot\Theta\left(m\right)\cdot h\label{trgm}
\end{equation}
for a unique $h\in G_{s\left(\pi\left(m\right)\right)}$. And, in
infinitesimal terms around $g=e$, 
\[
\Theta_{*,m}\circ\left(\rho_{m}\right)_{*,e}\left(\eta\right)=\left(R_{\Theta\left(m\right)}\right)_{*,e}\left(\eta\right)+\left(L_{\Theta\left(m\right)}\right)_{*,e}\left(\xi\right),
\]
for some $\xi\in\mathfrak{g}_{s\left(\pi\left(m\right)\right)}$.
In particular, if $m=s\left(\lambda\right)$ (see \eqref{tsl}), 
\begin{equation}
\Theta_{*,s\left(\lambda\right)}\circ\left(\rho_{s\left(\lambda\right)}\right)_{*,e}\left(\eta\right)=\eta+\xi.\label{trgm2}
\end{equation}

\bigskip{}

As we anticipate at the end of the last subsection, the submersions
$\Theta$ above defined play a role similar to that of $\psi$ in
a trivialization of a principal bundle (see \eqref{trivfi}). This
follows, for instance, by comparing \eqref{rots0} and \eqref{rots}.
In particular, we can see the map 
\[
A:v\in T_{m}U\longmapsto\left(R_{\Theta\left(m\right)}\right)_{*,e}^{-1}\Theta_{*,m}\left(v\right)\in\mathfrak{g}
\]
as some kind of flat connection for the submersion $\pi_{\left|U\right.}$.
Nevertheless, \eqref{c1} and \eqref{c2} are not satisfied in general.
In fact, we have from \eqref{trgm} that (for $g\in G$ and $m\in U$
such that $\rho\left(g,m\right)\in U$) 
\[
A\left(\left(\rho_{m}\right)_{*,e}\left(\eta\right)\right)=\eta+Ad_{\Theta\left(m\right)}\xi,\;\;\;\textrm{for some}\;\;\;\xi\in\mathfrak{g}_{m},
\]
and 
\[
A\left(\left(\rho_{g}\right)_{*,m}\left(v\right)\right)=Ad_{g}\left(A\left(v\right)+Ad_{\Theta\left(m\right)}\xi\right),\;\;\;\textrm{for some}\;\;\;\xi\in\mathfrak{g}_{m}.
\]

\subsection{Vertical and horizontal vector fields}

Fix again a point $m_{0}\in\mathcal{R}_{\rho}$. \begin{defn} \label{verhor}
We shall say that $X\in\mathfrak{X}\left(M\right)$ is \textbf{vertical
around} $m_{0}$ if 
\[
X\left(m\right)\in\mathsf{Ker}\left(\pi_{\left|U\right.}\right)_{*,m},\;\;\;\forall m\in U,
\]
and that $X$ is \textbf{horizontal at} $m_{0}$ if 
\[
X\left(m_{0}\right)\notin\mathsf{Ker}\left(\pi_{\left|U\right.}\right)_{*,m_{0}},
\]
for some admissible neighborhood $U$ of $m_{0}$. Finally, we shall
say that $X$ is $\Theta$-\textbf{horizontal} if $\mathsf{Im}X_{\left|U\right.}\subseteq\mathsf{Ker}\Theta_{*}$
for some horizontal submersion $\Theta:U\rightarrow\Theta\left(U\right)$
and some admissible neighborhood $U$ of $m_{0}$. \end{defn}

From \eqref{kerim}, it is clear that if there exists a function $\eta:U\to\mathfrak{g}$
such that $X(m)=(\rho_{m})_{*e}(\eta(m)),$ for all $m\in U,$ then
$X$ is vertical along $U.$ We are interested in vertical fields
which are in addition $G$-invariant (see \eqref{ginv}). For them,
we have the next result. \begin{prop} Consider $X\in\mathfrak{X}\left(M\right)$
such that, for some function $\eta:M\rightarrow\mathfrak{g}$, 
\begin{equation}
X\left(m\right)=\left(\rho_{m}\right)_{*,e}\left(\eta\left(m\right)\right),\;\;\;\forall m\in M.\label{vereta}
\end{equation}
Then $X$ is $G$-invariant if and only if 
\begin{equation}
\eta\left(\rho_{g}\left(m\right)\right)=Ad_{g}\left(\eta\left(m\right)\right)+\xi_{g,m},\label{adeta}
\end{equation}
for some $\xi_{g,m}\in\mathfrak{g}_{\rho_{g}\left(m\right)}$. We
shall say that $\eta$ is $Ad$-\textbf{equivariant} if $\xi_{g,m}=0$
for all $g,m$. \end{prop} \textit{Proof.} Since $\rho_{\rho_{g}\left(m\right)}=\rho_{m}\circ R_{g}$
and $\rho_{g}\circ\rho_{m}=\rho_{m}\circ L_{g}$, then 
\[
\begin{array}{lll}
\left(\rho_{g}\right)_{*,m}^{-1}\circ X\left(\rho_{g}\left(m\right)\right) & = & \left(\rho_{m}\right)_{*,e}\circ\left(L_{g^{-1}}\right)_{*,g}\circ\left(R_{g}\right)_{*,e}\left(\eta\left(\rho_{g}\left(m\right)\right)\right)\\
 & = & \left(\rho_{m}\right)_{*,e}\circ Ad_{g^{-1}}\left(\eta\left(\rho_{g}\left(m\right)\right)\right).
\end{array}
\]
Hence, \eqref{ginv} is fulfilled if and only if 
\[
Ad_{g^{-1}}\left(\eta\left(\rho_{g}\left(m\right)\right)\right)-\eta\left(m\right)\in\mathsf{Ker}\left(\left(\rho_{m}\right)_{*,e}\right)=\mathfrak{g}_{m},
\]
and the proposition follows from the fact that $Ad_{g}\left(\mathfrak{g}_{m}\right)=\mathfrak{g}_{\rho_{g}\left(m\right)}$.
$\;\;\;\blacksquare$

\noindent \bigskip{}

Regarding horizontal fields, note that if $X$ is $\Theta$-horizontal
and $X\left(m_{0}\right)\neq0$, then $X$ is horizontal at $m_{0}$.
Reciprocally, we have the next result. \begin{prop} \label{exhor}
If $X$ is horizontal at $m_{0}$ and $G$-invariant, then there exist
an admissible neighborhood $U$ of $m_{0}$, a section $s:\pi\left(U\right)\rightarrow U$
of $\pi_{\left|U\right.}$ satisfying \eqref{spi} and a horizontal
submersion $\Theta:U\rightarrow\Theta\left(U\right)$ such that $X$
is $\Theta$-horizontal\textbf{.} \end{prop} \textit{ Proof. } According
to Proposition 4.13 of \cite{gp}, if $X\left(m_{0}\right)\notin\mathsf{Ker}\left(\pi_{\left|U\right.}\right)_{*,m_{0}}$
for some admissible neighborhood $U$ of $m_{0}$, then, shrinking
$U$ if necessary, there exists a submersion $F:U\rightarrow F\left(U\right)$
transverse to $\pi_{\left|U\right.}$ such that 
\begin{equation}
X\left(m\right)\in\mathsf{Ker}F_{*,m},\;\;\;\forall m\in U.\label{xf}
\end{equation}
On the one hand, shrinking $U$ again, this gives rise to a diffeomorphism
\[
\mathfrak{D}\coloneqq\left(\pi_{\left|U\right.},F\right):U\rightarrow\pi\left(U\right)\times F\left(U\right).
\]
In terms of the latter, we have the section $s:\pi\left(U\right)\rightarrow U$
of $\pi_{\left|U\right.}$ given by 
\[
s\left(\pi\left(m\right)\right)=\mathfrak{D}^{-1}\left(\pi\left(m\right),F\left(m_{0}\right)\right),
\]
which satisfies $s\left(\pi\left(m_{0}\right)\right)=m_{0}$. So,
we have a section of $\pi_{\left|U\right.}$ satisfying \eqref{spi}
and, according to Theorem \ref{pro-tita}, this enable us to construct
a horizontal submersion $\Theta:U\rightarrow\Theta\left(U\right)$.
On the other hand, writing $s\left(\pi\left(m\right)\right)=\tilde{m}$,
\eqref{xf} says that 
\[
\mathfrak{D}_{*,\tilde{m}}\left(X\left(\tilde{m}\right)\right)=\left(\left(\pi_{\left|U\right.}\right)_{*,\tilde{m}}\left(X\left(\tilde{m}\right)\right),0\right),
\]
or equivalently 
\begin{equation}
X\left(\tilde{m}\right)=\left(\mathfrak{D}^{-1}\right)_{*,\left(\pi\left(m\right),F\left(m_{0}\right)\right)}\left(\left(\pi_{\left|U\right.}\right)_{*,\tilde{m}}\left(X\left(\tilde{m}\right)\right),0\right)=s_{*,\pi\left(m\right)}\circ\left(\pi_{\left|U\right.}\right)_{*,\tilde{m}}\left(X\left(\tilde{m}\right)\right).\label{xmono}
\end{equation}
In addition, the fact that $X$ is $G$-invariant ensures that (combine
\eqref{ginv} and \eqref{rots}) 
\[
\left(\rho_{\Theta\left(m\right)}\right)_{*,\tilde{m}}X\left(\tilde{m}\right)=X\left(m\right),
\]
and consequently (see \eqref{piro}) 
\[
\left(\pi_{\left|U\right.}\right)_{*,m}\left(X\left(m\right)\right)=\left(\pi_{\left|U\right.}\right)_{*,m}\circ\left(\rho_{\Theta\left(m\right)}\right)_{*,\tilde{m}}X\left(\tilde{m}\right)=\left(\pi_{\left|U\right.}\right)_{*,\tilde{m}}\left(X\left(\tilde{m}\right)\right)
\]
and (applying $\left(\rho_{\Theta\left(m\right)}\right)_{*,\tilde{m}}$
to \eqref{xmono}) 
\[
X\left(m\right)=\left(\rho_{\Theta\left(m\right)}\right)_{*,\tilde{m}}\circ s_{*,\pi\left(m\right)}\circ\left(\pi_{\left|U\right.}\right)_{*,m}\left(X\left(m\right)\right).
\]
Finally, using \eqref{vkt}, it follows that $\mathsf{Im}X_{\left|U\right.}\subseteq\mathsf{Ker}\Theta_{*}$,
as wanted. $\;\;\;\blacksquare$

\subsection{Local complete solutions from general group actions}

\noindent From above results and the duality between complete solutions
and first integrals, the theorem below easily follows.

\begin{thm}

\label{1}Fix $m_{0}\in\mathcal{R}_{\rho}$. 
\begin{enumerate}
\item If $X$ is vertical around $m_{0}$, then there exists an admissible
neighborhood $U$ of $m_{0}$ such that, for every section $s:\pi\left(U\right)\rightarrow U$
of $\pi_{\left|U\right.}$ satisfying \eqref{spi} and every $s$-horizontal
submersion $\Theta:U\rightarrow\Theta\left(U\right)$, the map 
\[
\Sigma\coloneqq\left(\Theta,\pi_{\left|U\right.}\right)^{-1}=\rho\circ\left(id_{\Theta\left(U\right)}\times s\right):\Theta\left(U\right)\times\pi\left(U\right)\rightarrow U
\]
is a complete solution of the $\Theta$-HJE for $X_{\left|U\right.}$. 
\item If $X$ is horizontal at $m_{0}$ and $G$-invariant, then there exist
an admissible neighborhood $U$ of $m_{0}$, a section $s:\pi\left(U\right)\rightarrow U$
of $\pi_{\left|U\right.}$ satisfying \eqref{spi} and a $s$-horizontal
submersion $\Theta:U\rightarrow\Theta\left(U\right)$ such that\footnote{By $\tau$ we are denoting the flipping map $\tau\left(x,y\right)=\left(y,x\right)$.}
\[
\Sigma\coloneqq\left(\pi_{\left|U\right.},\Theta\right)^{-1}=\rho\circ\tau\circ\left(s\times id_{\Theta\left(U\right)}\right):\pi\left(U\right)\times\Theta\left(U\right)\rightarrow U
\]
is a complete solution of the $\pi_{\left|U\right.}$-HJE for $X_{\left|U\right.}$. 
\end{enumerate}
\end{thm}

\noindent \textit{Proof. } In the first case we have that $\mathsf{Im}X_{\left|U\right.}\subseteq\mathsf{Ker}\left(\pi_{\left|U\right.}\right)_{*}$
and that $\pi_{\left|U\right.}$ and $\Theta$ are transverse. Using
the results of Section \ref{dual}, it follows that, shrinking $U$
if needed, $\Sigma\coloneqq\left(\Theta,\pi_{\left|U\right.}\right)^{-1}$
is a complete solution of the $\Theta$-HJE for $X_{\left|U\right.}$.
We only need to show that $\left(\Theta,\pi_{\left|U\right.}\right)^{-1}=\rho\circ\left(id_{\Theta\left(U\right)}\times s\right)$
. But, form \eqref{rots}, we have for all $m\in U$ that 
\[
\rho\circ\left(id_{\Theta\left(U\right)}\times s\right)\circ\left(\Theta,\pi_{\left|U\right.}\right)\left(m\right)=\rho\left(\Theta\left(m\right),s\left(\pi\left(m\right)\right)\right)=m.
\]

The second case can be proved in the same way, but using in addition
Proposition \ref{exhor} in order to ensure the existence of the section
$s$ and the submersion $\Theta$ such that $\mathsf{Im}X_{\left|U\right.}\in\mathsf{Ker}\Theta_{*}$.
$\;\;\;\blacksquare$

\begin{rem*} Regarding the objects described in Section \ref{tlf},
it is clear that the complete solutions $\Sigma$ given in the last
theorem, or more precisely their inverses $\Sigma^{-1}$, define the
analogue of a trivialization $\Psi:U\rightarrow\pi\left(U\right)\times G$
of a principal bundle. \end{rem*}

Summarizing, given a vertical vector field $X$ around $m_{0}\in\mathcal{R}_{\rho}$,
an admissible neighborhood $U$ of $m_{0}$ and a smooth section $s:\pi(U)\to U$
of $\pi_{|U},$ we have shown that a submersion $\Theta:U\rightarrow\Theta\left(U\right)$
and a complete solution of the $\Theta$-HJE for $X_{\left|U\right.}$
can be constructed up to quadratures. Also, given a horizontal vector
field $X$ at $m_{0}$, if $X$ is $G$-invariant, then there exists
a complete solution of the $\pi_{\left|U\right.}$-HJE for $X_{\left|U\right.}$.
But the latter has not been constructed up to quadratures (the proof
of Proposition 4.13 of \cite{gp}, which is used in Proposition \ref{exhor},
is based on the rectification of the field $X$).

\section{Horizontal dynamical systems and reconstruction}

Consider again a manifold $M$, a vector field $X\in\mathfrak{X}\left(M\right)$
and a group action $\rho:G\times M\rightarrow M$. Assume by now that
$\rho$ is free and proper, what implies that $\pi:M\rightarrow\left.M\right/G$
defines a principal fiber bundle. Assume also that $X$ is $G$-invariant,
and consequently $\pi$-related with a unique vector field $Y\in\mathfrak{X}\left(\left.M\right/G\right)$,
i.e. $\pi_{*}\circ X=Y\circ\pi$. In many cases, the integral curves
of $Y$ are known, and one is interested in constructing the integral
curves of $X$ from those of $Y$. Any procedure that enable us to
do that is usually called \textit{reconstruction}. The purpose of
this section is to show that there exists a deep connection between
reconstruction procedures and the complete solutions of a horizontal
vector field presented in Theorem \ref{1}, even though the action
$\rho$ is neither free nor proper.

\subsection{The usual reconstruction process}

Assume that we are in the setting of the beginning of this section
and we want to find the integral curve $\Gamma$ of $X$ such that
$\Gamma\left(0\right)=p_{0}$. Then we can (see, for instance, \cite{mont}): 
\begin{enumerate}
\item consider the integral curve $\gamma\left(t\right)$ of $Y$ such that
$\gamma\left(0\right)=\pi\left(p_{0}\right)$; 
\item fix a principal connection $A:TM\rightarrow\mathfrak{g}$; 
\item find a curve $d\left(t\right)$ such that 
\begin{equation}
A\left(d'\left(t\right)\right)=0\;\;\;\textrm{and}\;\;\;\pi\left(d\left(t\right)\right)=\gamma\left(t\right),\label{prerec}
\end{equation}
i.e. $d\left(t\right)$ is the horizontal lift of $\gamma\left(t\right)$; 
\item find $g\left(t\right)$ such that 
\begin{equation}
g'\left(t\right)=\left(L_{g\left(t\right)}\right)_{*,e}\left(\xi\left(t\right)\right)\;\;\;\textrm{and}\;\;\;g\left(0\right)=g_{0},\label{recec}
\end{equation}
with $\xi\left(t\right)=A\left(X\left(d\left(t\right)\right)\right)$
and $g_{0}$ such that $p_{0}=\rho\left(g_{0},d\left(0\right)\right)$. 
\end{enumerate}
It is easy to show that $\Gamma\left(t\right)=\rho\left(g\left(t\right),d\left(t\right)\right)$
is the integral curve we are looking for. The four steps above constitute
the usual \textbf{\textit{reconstruction process}}, and \eqref{prerec}
and \eqref{recec} the related\textbf{\textit{ reconstruction problem}}.\bigskip{}

If $X$ is vertical along all of $M$ (in the usual sense), i.e. $\mathsf{Im}X\subseteq\mathsf{Ker}\pi_{*}$,
then $Y=0$ and consequently the curves $d\left(t\right)$ and $\xi\left(t\right)$
are constant. In this case, we only have to solve \eqref{recec},
whose solutions are given by the exponential curves. We shall consider
this situation in the next section. So, suppose that $X\left(m\right)\notin\mathsf{Ker}\pi_{*,m}$,
for all $m\in M$. In that case we can consider a connection $A$
such that $X\in\mathsf{Ker}A$, i.e. $X$ is horizontal with respect
to $A$ (in the usual sense). Then, $\xi\left(t\right)=0$ and $g\left(t\right)=g_{0}$
for all $t$. Consequently, the reconstruction problem reduces to
solve \eqref{prerec}. In other words, we have the following alternative
(three steps) reconstruction process: 
\begin{enumerate}
\item consider the integral curve $\gamma\left(t\right)$ of $Y$ such that
$\gamma\left(0\right)=\pi\left(p_{0}\right)$; 
\item find a principal connection $A:TM\rightarrow\mathfrak{g}$ such that
$X$ is horizontal; 
\item find a curve $d\left(t\right)$ satisfying \eqref{prerec}. 
\end{enumerate}
Then, the curve $\Gamma\left(t\right)=\rho\left(g_{0},d\left(t\right)\right)$,
with $g_{0}$ such that $p_{0}=\rho\left(g_{0},d\left(0\right)\right)$,
is the integral curve of $X$ through $p_{0}$. In the next subsection,
we shall extend this procedure to Lie group actions which are not
necessarily free and proper.

\subsection{Reconstruction from complete solutions}

Let us go back to the general setting: a manifold $M$, a vector field
$X\in\mathfrak{X}\left(M\right)$ and a general Lie group action $\rho:G\times M\rightarrow M$.
Assume that $X$ is $G$-invariant and horizontal at every $m_{0}\in\mathcal{R}_{\rho}$
(see Definition \ref{verhor}). According to the second part of Theorem
\ref{1}, there exist an admissible neighborhood $U$ of $m_{0}$,
a section $s:\pi\left(U\right)\rightarrow U$ of $\pi_{\left|U\right.}$
satisfying \eqref{spi} and a $s$-horizontal submersion $\Theta:U\rightarrow\Theta\left(U\right)$
such that 
\[
\Sigma\coloneqq\rho\circ\tau\circ\left(s\times id_{\Theta\left(U\right)}\right):\pi\left(U\right)\times\Theta\left(U\right)\rightarrow U
\]
is a complete solution of the $\pi_{\left|U\right.}$-HJE for $X_{\left|U\right.}$.
The related partial solutions are functions 
\[
\sigma_{g}:\pi\left(U\right)\rightarrow U,\;\;\;g\in\Theta\left(U\right),
\]
such that $\sigma_{g}\left(\lambda\right)=\rho\left(g,s\left(\lambda\right)\right)$
for all $\lambda\in\pi\left(U\right)$ (see \eqref{psr}). In other
words, 
\begin{equation}
\sigma_{g}=\rho_{g}\circ s,\;\;\;g\in\Theta\left(U\right).\label{sgi}
\end{equation}

\begin{thm} Each vector field $X^{\sigma_{g}}\in\mathfrak{X}\left(\pi\left(U\right)\right)$
(see \eqref{Xs}) is equal, for all $g\in\Theta\left(U\right)$, to
the unique vector field $Y\in\mathfrak{X}\left(\pi\left(U\right)\right)$
such that 
\begin{equation}
\left(\pi_{\left|U\right.}\right)_{*}\circ X_{\left|U\right.}=Y\circ\pi_{\left|U\right.}.\label{pxyp}
\end{equation}
\end{thm} {\textit{Proof. } The Proposition \ref{pm} ensures the
existence of a unique vector field $Y\in\mathfrak{X}\left(\pi\left(U\right)\right)$
satisfying \eqref{pxyp}. So, we only must prove that $Y=X^{\sigma_{g}}$
for all $g\in\Theta\left(U\right)$. But from \eqref{Xs}, \eqref{sgi}
and \eqref{pxyp}, 
\[
\begin{array}{lll}
X^{\sigma_{g}} & = & \left(\pi_{\left|U\right.}\right)_{*}\circ X_{\left|U\right.}\circ\sigma_{g}=\left(\pi_{\left|U\right.}\right)_{*}\circ X_{\left|U\right.}\circ\rho_{g}\circ s\\
 & = & Y\circ\pi_{\left|U\right.}\circ\rho_{g}\circ s=Y\circ\pi_{\left|U\right.}\circ s=Y,
\end{array}
\]
as we wanted to show. $\;\;\;\blacksquare$

\noindent \bigskip{}

According to \eqref{G}, the integral curves $\Gamma$ of $X$ are
given by 
\[
\Gamma\left(t\right)=\sigma_{g}\left(\gamma\left(t\right)\right)=\rho\left(g,s\left(\gamma\left(t\right)\right)\right),
\]
where $\gamma$ is an integral curve of $Y=X^{\sigma_{g}}$. In other
words, above formula enable us to construct the integral curves of
$X$ from those of a vector field in the quotient. Note that $\pi\left(\Gamma\left(t\right)\right)=\gamma\left(t\right)$
and $\Theta\left(\Gamma\left(t\right)\right)=g$ for all $t$. Then,
given $p_{0}\in U$, in order to find the integral curve $\Gamma$
of $X_{\left|U\right.}$ such that $\Gamma\left(0\right)=p{}_{0}$,
we have the following (two steps) reconstruction process: 
\begin{enumerate}
\item consider the integral curve $\gamma\left(t\right)$ of $Y$ such that
$\gamma\left(0\right)=\pi\left(p_{0}\right)$; 
\item find a submersion $\Theta:U\rightarrow\Theta\left(U\right)$ such
that $X$ is $\Theta$-horizontal. 
\end{enumerate}
The curve 
\[
\Gamma\left(t\right)=\Sigma\left(\gamma\left(t\right),g_{0}\right)=\rho\left(g_{0},s\left(\gamma\left(t\right)\right)\right),\mbox{ with \ensuremath{g_{0}=\Theta\left(p_{0}\right)},}
\]
is the one we are looking for. So, the complete solution $\Sigma$
solves the reconstruction problem (around $m_{0}$).

\section{Vertical dynamical systems and integrability by quadratures}

In this section, using the integrability result of Section $2$ (see
Theorem \ref{0}), we show that the exponential curves $t\mapsto\exp\left(\xi\,t\right)$
of a Lie group $G$, for some points $\xi$ of its Lie algebra $\mathfrak{g}$,
can be explicitly constructed up to quadratures. Moreover, we show
that, for compact and for semisimple Lie groups, such a construction
works for all $\xi$ inside a dense open subset of $\mathfrak{g}$.
Then, we state sufficient conditions under which a vertical (and invariant)
vector field is integrable up to quadratures.

\subsection{Invariant and vertical vector fields}

\label{ivf} Consider again a manifold $M$, a vector field $X\in\mathfrak{X}\left(M\right)$
and a Lie group action $\rho:G\times M\rightarrow M$. Assume that
$X$ is vertical around every $\rho$-regular point $m_{0}$ (see
Definition \ref{verhor}) and consider a complete solution 
\begin{equation}
\Sigma\coloneqq\rho\circ\left(id_{\Theta\left(U\right)}\times s\right):\Theta\left(U\right)\times\pi\left(U\right)\rightarrow U\label{S}
\end{equation}
as those given in the first part of Theorem \ref{1}. The related
partial solutions are 
\[
\sigma_{\lambda}:\Theta\left(U\right)\subseteq G\rightarrow U,\;\;\;\lambda\in\pi\left(U\right),
\]
with $\sigma_{\lambda}\left(g\right)=\rho\left(g,s\left(\lambda\right)\right)$
for all $g\in\Theta\left(U\right)$ (see \eqref{psr}). In other words,
\[
\sigma_{\lambda}=\rho_{s\left(\lambda\right)},\;\;\;\lambda\in\pi\left(U\right).
\]

\begin{thm} \label{2}If $X$ is $G$-invariant, then the vector
field $X^{\sigma_{\lambda}}\in\mathfrak{X}\left(\Theta\left(U\right)\right)$
(see \eqref{Xs}) is given by 
\begin{equation}
X^{\sigma_{\lambda}}\left(g\right)=\left(L_{g}\right)_{*,e}\left(\eta_{\lambda}\right)\label{Xsle}
\end{equation}
for a unique vector (recall \eqref{tsl}) 
\begin{equation}
\eta_{\lambda}=\Theta_{*,s\left(\lambda\right)}\circ X\left(s\left(\lambda\right)\right)\in T_{e}\Theta\left(U\right)\subseteq\mathfrak{g}.\label{el}
\end{equation}
In particular, if $X$ is given by \eqref{vereta} and \eqref{adeta},
then 
\begin{equation}
\eta_{\lambda}=\eta\left(s\left(\lambda\right)\right)+\xi_{\lambda}\label{eles}
\end{equation}
for some $\xi_{\lambda}\in\mathfrak{g}_{s\left(\lambda\right)}$.
In any case, the integral curve of $X$ passing through $p_{0}\in U$
at $t=0$ can be written 
\begin{equation}
\Gamma\left(t\right)=\rho\left(g_{0}\,\exp\left(\left(\eta_{\lambda}+\chi_{\lambda}\right)\,t\right),s\left(\lambda\right)\right),\label{Gtrg}
\end{equation}
with $\left(g_{0},\lambda\right)=\left(\Theta\left(p_{0}\right),\pi\left(p_{0}\right)\right)$
and $\chi_{\lambda}\in\mathfrak{g}_{s\left(\lambda\right)}$ arbitrary.
\end{thm}

\noindent \textit{Proof. } We know that (see \eqref{hjrel}) 
\begin{equation}
X\circ\sigma_{\lambda}=\left(\sigma_{\lambda}\right)_{*}\circ X^{\sigma_{\lambda}},\label{srel}
\end{equation}
and consequently $\mathsf{Im}\left(X\circ\sigma_{\lambda}\right)\subseteq\mathsf{Im}\left(\sigma_{\lambda}\right)_{*}$.
Since $\sigma_{\lambda}$ is a diffeomorphism onto its image, then
\[
\left(\sigma_{\lambda}\right)_{*,e}:T_{e}\Theta\left(U\right)\rightarrow T_{s\left(\lambda\right)}\mathsf{Im}\sigma_{\lambda}=\mathsf{Im}\left(\sigma_{\lambda}\right)_{*,e}
\]
is a linear isomorphism. So, for $X\circ\sigma_{\lambda}\left(e\right)=X\left(s\left(\lambda\right)\right)\in\mathsf{Im}\left(\sigma_{\lambda}\right)_{*,e}$
there exists a unique vector $\eta_{\lambda}\in\mathfrak{g}$ such
that 
\[
\left(\sigma_{\lambda}\right)_{*,e}\left(\eta_{\lambda}\right)=X\left(s\left(\lambda\right)\right).
\]
Let us apply $\left(\rho_{g}\right)_{*,s\left(\lambda\right)}$ on
both members of above equation. For the first member we have that
\[
\left(\rho_{g}\right)_{*,s\left(\lambda\right)}\circ\left(\sigma_{\lambda}\right)_{*,e}\left(\eta_{\lambda}\right)=\left(\sigma_{\lambda}\right)_{*,g}\circ\left(L_{g}\right)_{*,e}\left(\eta_{\lambda}\right),
\]
where we have used that $\sigma_{\lambda}=\rho_{s\left(\lambda\right)}$
and the identity 
\[
\rho_{g}\circ\rho_{s\left(\lambda\right)}=\rho_{s\left(\lambda\right)}\circ L_{g}.
\]
For the second member, using the $G$-invariance of $X$ (recall \eqref{ginv}),
we have that 
\[
\left(\rho_{g}\right)_{*,s\left(\lambda\right)}\left(X\left(s\left(\lambda\right)\right)\right)=X\left(\rho_{g}\left(s\left(\lambda\right)\right)\right)=X\left(\sigma_{\lambda}\left(g\right)\right).
\]
Then 
\[
\left(\sigma_{\lambda}\right)_{*,g}\circ\left(L_{g}\right)_{*,e}\left(\eta_{\lambda}\right)=X\left(\sigma_{\lambda}\left(g\right)\right)
\]
and, consequently, \eqref{Xsle} follows from \eqref{srel} and the
injectivity of $\sigma_{\lambda}$. Finally, using \eqref{Xsle} and
the fact that $X^{\sigma_{\lambda}}=\Theta_{*}\circ X\circ\sigma_{\lambda}$
(see \eqref{Xs}), 
\[
\eta_{\lambda}=X^{\sigma_{\lambda}}\left(e\right)=\Theta_{*}\circ X\circ\sigma_{\lambda}\left(e\right)=\Theta_{*,s\left(\lambda\right)}\circ X\left(s\left(\lambda\right)\right),
\]
which gives precisely \eqref{el}. In particular, if $X$ is given
by \eqref{vereta} and \eqref{adeta}, using \eqref{trgm2} we obtain
easily \eqref{eles}.\bigskip{}

Now, let us prove \eqref{Gtrg}. Given a curve 
\[
\Gamma\left(t\right)=\rho\left(\gamma\left(t\right),s\left(\lambda\right)\right)=\rho_{s\left(\lambda\right)}\left(\gamma\left(t\right)\right)=\sigma_{\lambda}\left(\gamma\left(t\right)\right),
\]
with $\gamma\left(t\right)=g_{0}\,\exp\left(\left(\eta_{\lambda}+\chi_{\lambda}\right)\,t\right)$,
since $\gamma'\left(t\right)=\left(L_{\gamma\left(t\right)}\right)_{*,e}\left(\eta_{\lambda}+\chi_{\lambda}\right)$,
then 
\[
\Gamma'\left(t\right)=\left(\sigma_{\lambda}\right)_{*,\gamma\left(t\right)}\circ\left(L_{\gamma\left(t\right)}\right)_{*,e}\left(\eta_{\lambda}\right)+\left(\sigma_{\lambda}\right)_{*,\gamma\left(t\right)}\circ\left(L_{\gamma\left(t\right)}\right)_{*,e}\left(\chi_{\lambda}\right).
\]
On the other hand, given a curve $h\left(x\right)$ on $G_{s\left(\lambda\right)}$
such that $h\left(0\right)=e$ and $h'\left(0\right)=\chi_{\lambda}$,
since 
\[
\sigma_{\lambda}\left(L_{g}\left(h\left(x\right)\right)\right)=\rho\left(g\,h\left(x\right),s\left(\lambda\right)\right)=\rho\left(g,s\left(\lambda\right)\right),\;\;\;\forall g\in G,
\]
then $\left(\sigma_{\lambda}\right)_{*,\gamma\left(t\right)}\circ\left(L_{\gamma\left(t\right)}\right)_{*,e}\left(\chi_{\lambda}\right)=0$
for all $t$. Accordingly, using \eqref{Xsle} and \eqref{srel},
\[
\Gamma'\left(t\right)=\left(\sigma_{\lambda}\right)_{*,\gamma\left(t\right)}\circ\left(L_{\gamma\left(t\right)}\right)_{*,e}\left(\eta_{\lambda}\right)=\left(\sigma_{\lambda}\right)_{*,\gamma\left(t\right)}\circ X^{\sigma_{\lambda}}\left(\gamma\left(t\right)\right)=X\left(\sigma_{\lambda}\left(\gamma\left(t\right)\right)\right)=X\left(\Gamma\left(t\right)\right),
\]
i.e. $\Gamma$ is an integral curve of $X$. If $\Gamma\left(0\right)=p_{0}$,
then (see \eqref{S}) 
\[
p_{0}=\sigma_{\lambda}\left(g_{0}\right)=\Sigma\left(g_{0},\lambda\right)=\left(\Theta,\pi_{\left|U\right.}\right)^{-1}\left(g_{0},\lambda\right),
\]
what ends our proof. $\;\;\;\blacksquare$

\bigskip{}

In order to consider concrete examples of vertical and $G$-invariant
fields, suppose that $M$ is a symplectic manifold, with symplectic
form $\omega$, and $\rho$ is a symplectic action with an $Ad^{*}$-equivariant
momentum map $K$.

\begin{prop} \label{7} If $\phi:\mathfrak{g}^{*}\rightarrow\mathfrak{g}$
is \textbf{equivariant}, i.e. 
\begin{equation}
Ad_{g}\phi\left(\alpha\right)=\phi\left(Ad_{g}^{*}\alpha\right),\;\;\;\forall\alpha\in\mathfrak{g}^{*},\;g\in G,\label{adfi}
\end{equation}
then 
\[
X=\omega^{\sharp}\circ K^{*}\phi
\]
is vertical and $G$-invariant. Also, given an admissible neighborhood
$U$ and a complete solution \eqref{S}, the related vector $\eta_{\lambda}\in\mathfrak{g}$
(see \eqref{el}), for each $\lambda\in\pi\left(U\right)$, is given
by 
\begin{equation}
\eta_{\lambda}=\phi\left(K\left(s\left(\lambda\right)\right)\right)+\xi_{\lambda},\label{elfi}
\end{equation}
for some $\xi_{\lambda}\in\mathfrak{g}_{s\left(\lambda\right)}$.
In particular, the integral curve of $X$ passing through $p_{0}\in U$
at $t=0$ can be written 
\begin{equation}
\Gamma\left(t\right)=\rho\left(g_{0}\,\exp\left(\phi\left(K\left(s\left(\lambda\right)\right)\right)\,t\right),s\left(\lambda\right)\right),\label{Gtrg-1}
\end{equation}
with $\left(g_{0},\lambda\right)=\left(\Theta\left(p_{0}\right),\pi\left(p_{0}\right)\right)$.
\end{prop}

\noindent \textit{Proof. }The form of $X$ ensures that $X\left(m\right)\in\left(\mathsf{Ker}K_{*,m}\right)^{\bot}$,
for all $m\in M$. Then, for every admissible neighborhood $U$, we
have from \eqref{kkkp} that 
\[
X\left(m\right)\in\left(\mathsf{Ker}\left(\pi_{\left|U\right.}\right)_{*,m}\right),\quad\forall m\in U,
\]
i.e. $X$ is vertical. To show $G$-invariance, note first that \eqref{adinv}
implies the equality 
\[
K_{*,\rho_{g}\left(m\right)}\circ\left(\rho_{g}\right)_{*,m}=Ad_{g}^{*}\circ K_{*,m},
\]
and in dual form (changing $g$ by $g^{-1}$ and $m$ by $\rho_{g}\left(m\right)$)
\[
\left(\rho_{g^{-1}}\right)_{\rho_{g}\left(m\right)}^{*}\circ K_{m}^{*}=K_{\rho_{g}\left(m\right)}^{*}\circ Ad_{g}.
\]
Combining the equation above, \eqref{sya}, \eqref{adinv} and \eqref{adfi},
we obtain 
\[
\left(\rho_{g}\right)_{*,m}\left(X\left(m\right)\right)=\left(\rho_{g}\right)_{*,m}\left(\omega_{m}^{\sharp}\circ K_{m}^{*}\left(\phi\left(K\left(m\right)\right)\right)\right)=X\left(\rho_{g}\left(m\right)\right),
\]
as desired. Now, we will shown \eqref{elfi}. From the very definition
of $K$ (see \eqref{defk}) we have that $\omega_{m}^{\flat}\circ\left(\rho_{m}\right)_{*,e}=K_{m}^{*}$,
so 
\[
\eta_{\lambda}=\Theta_{*,s\left(\lambda\right)}\circ X\left(s\left(\lambda\right)\right)=\Theta_{*,s\left(\lambda\right)}\circ\left(\rho_{s\left(\lambda\right)}\right)_{*,e}\left(\phi\left(K\left(s\left(\lambda\right)\right)\right)\right).
\]
Thus, \eqref{elfi} follows from \eqref{trgm2}. Finally, \eqref{Gtrg-1}
is a direct consequence of the previous theorem. $\;\;\;\blacksquare$

\subsection{The cotangent bundle and the left multiplication}

\subsubsection{A class of invariant vertical vectors}

\noindent Given a Lie group $G$, consider its cotangent bundle $T^{*}G$
with its canonical symplectic structure $\omega_{G}=-d\theta_{G}$.
Consider also the action 
\[
\rho:G\times T^{*}G\rightarrow T^{*}G
\]
such that, for all $g\in G$ and $\alpha_{h}\in T_{h}^{*}G$, 
\[
\rho\left(g,\alpha_{h}\right)=\left[\left(L_{g}\right)_{h}^{*}\right]^{-1}\left(\alpha_{h}\right)\in T_{gh}^{*}G.
\]
Note that $\rho$ is symplectic (see \eqref{sya}) and has an $Ad^{*}$-equivariant
momentum map $J:T^{*}G\rightarrow\mathfrak{g}^{*}$ given by 
\[
J\left(\alpha_{g}\right)=\left(R_{g}\right)_{e}^{*}\left(\alpha_{g}\right).
\]
Also, $\rho$ is a free and proper action, the quotient $\left.T^{*}G\right/G$
is a manifold diffeomorphic to $\mathfrak{g}^{*}$ and the canonical
projection $\pi$ can be seen as the submersion 
\[
\pi:T^{*}G\rightarrow\mathfrak{g}^{*}:\alpha_{g}\longmapsto\left(L_{g}\right)_{e}^{*}\left(\alpha_{g}\right).
\]
In other words, every point of $T^{*}G$ is $\rho$-regular and the
whole of $T^{*}G$ is an admissible neighborhood. Then, according
to Proposition \ref{pm} (see Eq. \eqref{T}), 
\begin{equation}
\Xi^{\sharp}\left(\eta\right)\coloneqq\pi_{*,\alpha}\circ\omega_{G}^{\sharp}\circ\pi_{\alpha}^{*}\left(\eta\right),\;\;\;\eta\in T_{\alpha}^{*}\mathfrak{g}^{*},\label{expc}
\end{equation}
defines a Poisson bracket on $\mathfrak{g}^{*}$ and $\pi$ is a Poisson
morphism between $\left(T^{*}G,\omega_{G}\right)$ (with its related
Poisson structure) and $\left(\mathfrak{g}^{*},\Xi\right)$. Moreover,
it can be shown that $\Xi$ is the \textbf{Kirillov-Kostant bracket}
on $\mathfrak{g}^{*}$ (see \cite{mr}), i.e. 
\begin{equation}
\Xi^{\sharp}\left(\eta\right)=ad_{\eta}^{*}\alpha,\;\;\;\eta\in T_{\alpha}^{*}\mathfrak{g}^{*}\cong\mathfrak{g}.\label{pb}
\end{equation}
On the other hand, the map $s:\mathfrak{g}^{*}\rightarrow T^{*}G$
such that $s\left(\alpha\right)=\alpha\in T_{e}^{*}G=\mathfrak{g}^{*}$
is a global section of $\pi$ and satisfies $s\left(\pi\left(\alpha\right)\right)=\alpha$
for all $\alpha\in\mathfrak{g}^{*}$. A related horizontal submersion
is the map $\Theta:T^{*}G\rightarrow G$ such that $\Theta\left(\alpha_{g}\right)=g$,
i.e. the canonical cotangent projection $\pi_{G}:T^{*}G\rightarrow G$.
In fact, 
\[
\rho\left(\pi_{G}\left(\alpha_{g}\right),s\left(\pi\left(\alpha_{g}\right)\right)\right)=\rho\left(g,\left(L_{g}\right)_{e}^{*}\left(\alpha_{g}\right)\right)=\left[\left(L_{g}\right)_{e}^{*}\right]^{-1}\left(\left(L_{g}\right)_{e}^{*}\left(\alpha_{g}\right)\right)=\alpha_{g},
\]
for all $\alpha_{g}\in T_{g}^{*}G$.

\begin{rem}\label{rm2}{ Note that $\left(\pi_{G},\pi\right):T^{*}G\rightarrow G\times\mathfrak{g}^{*}$
is the left trivialization of $T^{*}G.$ Thus, $T^{*}G$ may be identified
with $G\times\mathfrak{g}$ and, under this identification, the projections
$\pi_{G}:T^{*}G\to G$ and $\pi:T^{*}G\to\mathfrak{g}^{*}$ are just
the canonical projections 
\[
pr_{1}:G\times\mathfrak{g}^{*}\to G\quad\mbox{and}\quad pr_{2}:G\times\mathfrak{g}^{*}\to\mathfrak{g}^{*}
\]
on the first and second factor, respectively. Moreover, the canonical
symplectic structure $\omega_{G}$ on $T^{*}G$ is the $2$-form on
$G\times\mathfrak{g}^{*}$ given by 
\begin{equation}
\omega_{G}(g,\alpha)((v_{g},\beta),(v_{g}',\beta'))=\beta'((L_{g^{-1}})_{*g}(v_{g}))-\beta((L_{g^{-1}})_{*g}(v_{g}'))+(ad_{(L_{g^{-1}})_{*g}(v_{g})}^{*}\alpha)((L_{g^{-1}})_{*g}(v_{g}')),\label{89'}
\end{equation}
for $(g,\alpha)\in G\times\mathfrak{g}^{*}$ and $(v_{g},\beta),(v_{g}',\beta')\in T_{g}G\times\mathfrak{g}^{*}\cong T_{g}G\times T_{\alpha}\mathfrak{g}^{*}$
(see \cite{am}). In addition, the action $\rho:G\times(G\times\mathfrak{g}^{*})\to G\times\mathfrak{g}^{*}$
is just the left translation on the first factor, that is, 
\[
\rho\left(g,\left(g',\alpha\right)\right)=\left(gg',\alpha\right)\quad\mbox{for \ensuremath{g,g'\in G}\quad and\quad\ensuremath{\alpha\in\mathfrak{g}^{*}}},
\]
and the momentum map $J:G\times\mathfrak{g}\to\mathfrak{g}^{*}$ is
just the co-adjoint action of $G$ on $\mathfrak{g}^{*}$ 
\[
J\left(g,\alpha\right)=Ad_{g^{-1}}^{*}\alpha
\]
(for more details, see \cite{am}).} \end{rem}

According to Theorem \ref{1} (part 1), for every vertical vector
field $X\in\mathfrak{X}\left(T^{*}G\right)$ along all of $T^{*}G$,
\begin{equation}
\Sigma=\left(\pi_{G},\pi\right)^{-1}=\rho\circ\left(id_{G}\times s\right):G\times\mathfrak{g}^{*}\rightarrow T^{*}G\label{sG}
\end{equation}
is a (global) complete solution of the $\pi_{G}$-HJE for $X$. If
in addition $X$ is $G$-invariant, then its integral curve with initial
condition $\Gamma\left(0\right)=\left(\pi_{G},\pi\right)^{-1}\left(g_{0},\alpha\right)$
is (recall \eqref{el} and \eqref{Gtrg}) 
\[
\Gamma\left(t\right)=\rho\left(g_{0}\,\exp\left(\eta_{\alpha}\,t\right),s\left(\alpha\right)\right)=\left[\left(L_{g_{0}\,\exp\left(\eta_{\alpha}\,t\right)}\right)_{e}^{*}\right]^{-1}\left(\alpha\right),
\]
with 
\begin{equation}
\eta_{\alpha}=\left(\pi_{G}\right)_{*,\alpha}\left(X\left(\alpha\right)\right).\label{90'}
\end{equation}
We have used above that $s\left(\alpha\right)=\alpha$. That is the
case, for instance, of a vector field $X$ of the form 
\begin{equation}
X^{\phi}=\omega_{G}^{\sharp}\circ\pi^{*}\phi,\label{xfi}
\end{equation}
with $\phi\in\Omega^{1}\left(\mathfrak{g}^{*}\right)$ such that 
\begin{equation}
\Xi^{\sharp}\left(\phi\left(\alpha\right)\right)=ad_{\phi\left(\alpha\right)}^{*}\alpha=0,\label{cafi}
\end{equation}
i.e. $\phi$ is a Casimir $1$-form with respect to the Poisson bracket
\eqref{pb}. The $G$-invariance of $X^{\phi}$ is immediate from
\eqref{sya} and \eqref{ginv}, and the verticality is ensured by
\eqref{expc} and \eqref{cafi}. On the other hand, it can be shown
that $\eta_{\alpha}=\phi\left(\alpha\right)$. In fact, using the
left trivialization of $T^{*}G,$ we can identify $T^{*}G$ with $G\times\mathfrak{g}^{*}$
(see Remark \ref{rm2}). Under this identification, $X^{\phi}$ may
be considered as a vector field on $G\times\mathfrak{g}^{*}$ and,
from (\ref{89'}) and (\ref{xfi}), it follows that 
\[
X^{\phi}\left(g,\alpha\right)=\left(\left(L_{g}\right){}_{*,e}\left(\phi\left(\alpha\right)\right),ad_{\phi\left(\alpha\right)}^{*}\alpha\right)=\left(\left(L_{g}\right){}_{*,e}\left(\phi\left(\alpha\right)\right),0\right),
\]
for all $(g,\alpha)\in G\times{\frak{g}}^{*}$. Thus, using (\ref{90'})
and Remark \ref{rm2}, we deduce that

\[
\eta_{\alpha}=\left(\pi_{G}\right){}_{*,\alpha}\left(X_{\left(e,\alpha\right)}^{\phi}\right)=\phi\left(\alpha\right),\mbox{ for all }\alpha\in\mathfrak{g}^{*}.
\]

So, in terms of $\Sigma$, the trajectories of $X^{\phi}$ can be
written 
\[
\Gamma\left(t\right)=\Sigma\left(g_{0}\,\exp\left(\phi\left(\alpha\right)\,t\right),\alpha\right)=\rho\left(g_{0}\exp\left(\phi\left(\alpha\right)\,t),\alpha\right)\right).
\]

\subsubsection{Construction of the exponential curves up to quadratures}

In this subsection, we are going to show that $X^{\phi}$ (see \eqref{xfi})
is integrable up to quadratures (on a dense subset of $T^{*}G$) and,
consequently, the exponential curves $\exp\left(\phi\left(\alpha\right)\,t\right)$
can be explicitly obtained, also up to quadratures. The proof will
be based on Theorem \ref{0}.

\begin{prop} Consider the co-adjoint action $Ad^{*}:G\times\mathfrak{g}^{*}\rightarrow\mathfrak{g}^{*}$
and the related isotropy subgroups $G_{\alpha}$, with $\alpha\in\mathfrak{g}^{*}$.
Then, for every $Ad^{*}$-regular point $\alpha_{0}$ and any admissible
neighborhood $V\subseteq\mathfrak{g}^{*}$ of $\alpha_{0}$, the function
\begin{equation}
F\coloneqq\left(J,\pi\right)_{\left|U\right.}:U=J^{-1}(V)\subseteq T^{*}G\rightarrow\mathfrak{g}^{*}\times\mathfrak{g}^{*}\label{Fu}
\end{equation}
is a submersion onto the closed submanifold 
\[
F\left(U\right)\coloneqq\left\{ \left(\alpha,Ad_{g}^{*}\alpha\right)\;:\;\alpha\in V,\;g\in G\right\} \subseteq\mathfrak{g}^{*}\times\mathfrak{g}^{*}.
\]
Moreover, 
\begin{equation}
\mathsf{Ker}F_{*}=\mathsf{Ker}\left(J_{\left|U\right.}\right)_{*}\cap\mathsf{Ker}\left(\pi_{\left|U\right.}\right)_{*}\label{intk}
\end{equation}
and 
\begin{equation}
\mathsf{Ker}F_{*}\subseteq\left(\mathsf{Ker}F_{*}\right)^{\bot}.\label{isotf}
\end{equation}
\end{prop}

\noindent \textit{Proof.} It is easy to see that the composition of
$\left(J,\pi\right):T^{*}G\rightarrow\mathfrak{g}^{*}\times\mathfrak{g}^{*}$
and 
\[
\mathfrak{R}^{-1}:G\times\mathfrak{g}^{*}\rightarrow T^{*}G:\left(g,\alpha\right)\rightarrow\left(\left(R_{g}\right)_{e}^{*}\right)^{-1}\left(\alpha\right)
\]
(the inverse of the right trivialization) gives 
\[
\left(J,\pi\right)\circ\mathfrak{R}^{-1}\left(g,\alpha\right)=\left(\alpha,Ad_{g}^{*}\alpha\right).
\]
Then, given an $Ad^{*}$-regular point $\alpha_{0}$ and an admissible
neighborhood $V\subseteq\mathfrak{g}^{*}$ of $\alpha_{0}$, we have
from Proposition \ref{prop5'} (applied to the action $Ad^{*}$) that
$\left(J,\pi\right)\circ\mathfrak{R}^{-1}$ restricted to $G\times V$
is a submersion onto the closed submanifold $F\left(U\right)\subseteq V\times\mathfrak{g}^{*}$.
As a consequence, since $\mathfrak{R}^{-1}$ is a diffeomorphism,
the first affirmation of the proposition follows. On the other hand,
\eqref{intk} follows straightforwardly and \eqref{isotf} is a direct
consequence of the identity $\mathsf{Ker}J_{*}=\left(\mathsf{Ker}\pi_{*}\right)^{\bot}$
and the inclusion 
\[
\begin{array}{lll}
\left(\mathsf{Ker}J_{*}\cap\mathsf{Ker}\pi_{*}\right)^{\bot} & = & \left(\mathsf{Ker}J_{*}\right)^{\bot}+\left(\mathsf{Ker}\pi_{*}\right)^{\bot}\\
 & = & \mathsf{Ker}\pi_{*}+\mathsf{Ker}J_{*}\supseteq\mathsf{Ker}J_{*}\cap\mathsf{Ker}\pi_{*}.\;\;\;\blacksquare
\end{array}
\]

\bigskip{}

\begin{rem*} A similar result was proved in \cite{j} (see Theorem
4.1 there), but in terms of Poisson sub-algebras (see Remark \ref{ncij}).
\end{rem*}

Because of the form of $X^{\phi}$, it is clear that $\mathsf{Im}X^{\phi}\subseteq\left(\mathsf{Ker}\pi_{*}\right)^{\bot}$.
As a consequence (recall \eqref{kkkp}) 
\begin{equation}
\mathsf{Im}X^{\phi}\subseteq\mathsf{Ker}J_{*}\cap\mathsf{Ker}\pi_{*}.\label{imxpj}
\end{equation}
So, using the last proposition and combining \eqref{intk}, \eqref{isotf}
and \eqref{imxpj}, it follows that, for each $Ad^{*}$-regular point
$\alpha_{0}$, we can construct a neighborhood $U$ of $\alpha_{0}$
and a submersion $F:U\rightarrow F\left(U\right)$ (given by \eqref{Fu})
such that 
\[
\mathsf{Im}X_{\left|U\right.}^{\phi}\subseteq\mathsf{Ker}F_{*}\;\;\;\textrm{and}\;\;\;\mathsf{Ker}F_{*}\subseteq\left(\mathsf{Ker}F_{*}\right)^{\bot}.
\]

\begin{rem} \label{ncir} It can be shown that $\left(\mathsf{Ker}F_{*}\right)^{\bot}$
is an integrable distribution. Then, if $\phi$ is an exact $1$-form,
$X_{\left|U\right.}^{\phi}$ and $F$ define a NCI system on $U$
(see Section \ref{ncis}). \end{rem}

In addition, since $\mathsf{Ker}F_{*}\subseteq\mathsf{Ker}\left(\pi_{\left|U\right.}\right)_{*}$,
we have that 
\[
L_{X^{\phi}}\beta^{\phi}=0,\;\;\;\textrm{with}\;\;\;\beta^{\phi}=\pi^{*}\phi,
\]
as we saw at the end of Section \ref{intq} (recall \eqref{FG} and
\eqref{p3}). This enable us to apply Theorem \ref{0} to $X_{\left|U\right.}^{\phi}$.

\begin{rem*}For each $n\in\Pi\left(U\right)$, the map $\beta_{n}^{\phi}\in\Omega^{1}\left(F\left(U\right)\right)$
related to $\beta^{\phi}$ and defined by \eqref{dbn}, is given by
\begin{equation}
\beta_{n}^{\phi}=pr_{1}^{*}\phi,\label{bfg}
\end{equation}
with $pr_{1}:\mathfrak{g}^{*}\times\mathfrak{g}^{*}\rightarrow\mathfrak{g}^{*}$
the projection onto the first factor. \end{rem*}

According to the proof of Theorem \ref{0}, we can construct up to
quadratures a submersion $\Pi:U\rightarrow\Pi\left(U\right)$ transverse
to $F$, a complete solution $\hat{\Sigma}\coloneqq\left(\Pi,F\right)^{-1}$
and a family of immersions (see \eqref{fit}) 
\[
\varphi_{\lambda}:\Pi\left(U\right)\rightarrow T_{\lambda}^{*}F\left(U\right),\;\;\;\lambda\in F\left(U\right),
\]
such that the integral curves of $X_{\left|U\right.}^{\phi}$ are
given by $\Gamma\left(t\right)=\hat{\Sigma}\left(\gamma\left(t\right),\lambda\right)$,
with $\gamma$ satisfying (see \eqref{eg} and \eqref{bfg}) 
\begin{equation}
\varphi_{\lambda}\left(\gamma\left(t\right)\right)=\varphi_{\lambda}\left(\gamma\left(0\right)\right)+t\,\phi\left(pr_{1}\left(\lambda\right)\right).\label{0bfg}
\end{equation}
On the other hand, we know that $\Gamma\left(t\right)=\Sigma\left(g_{0}\,\exp\left(\phi\left(\alpha\right)\,t\right),\alpha\right)$
for some $g_{0}$ and $\alpha$, with $\Sigma$ given by \eqref{sG}.
So, for an integral curve passing through 
\begin{equation}
\Gamma\left(0\right)=\Sigma\left(e,\alpha\right)=\alpha=\hat{\Sigma}\left(\gamma\left(0\right),\lambda\right),\label{ic}
\end{equation}
since 
\[
\lambda=F\circ\hat{\Sigma}\left(\gamma\left(0\right),\lambda\right)=F\left(\alpha\right)=\left(J,\pi\right)\left(\alpha\right)=\left(\alpha,\alpha\right),
\]
we have that $\Sigma\left(\exp\left(\phi\left(\alpha\right)\,t\right),\alpha\right)=\hat{\Sigma}\left(\gamma\left(t\right),\left(\alpha,\alpha\right)\right)$.
Consequently 
\[
\exp\left(\phi\left(\alpha\right)\,t\right)=\pi_{G}\circ\hat{\Sigma}\left(\gamma\left(t\right),\left(\alpha,\alpha\right)\right),
\]
with $\gamma$ satisfying (see \eqref{0bfg}) 
\[
\varphi_{\left(\alpha,\alpha\right)}\left(\gamma\left(t\right)\right)=\varphi_{\left(\alpha,\alpha\right)}\left(\Pi\left(\alpha\right)\right)+t\,\phi\left(\alpha\right),
\]
and where we have used that (see \eqref{ic}) 
\[
\gamma\left(0\right)=\Pi\circ\hat{\Sigma}\left(\gamma\left(0\right),\lambda\right)=\Pi\circ\Sigma\left(e,\alpha\right)=\Pi\left(\alpha\right).
\]
Hence, we have shown the next result.

\begin{prop} \label{p6} Given a Casimir $1$-form $\phi:\mathfrak{g}^{*}\rightarrow\mathfrak{g}$
and a point $\alpha\in\mathcal{R}_{Ad^{*}}$, the exponential curve
$\exp\left(\phi\left(\alpha\right)\,t\right)$ can be constructed
up to quadratures. More explicitly, it is given by the formula 
\begin{equation}
\exp\left(\phi\left(\alpha\right)\,t\right)=\pi_{G}\left(\left(\Pi,F\right)^{-1}\left(\varphi_{\left(\alpha,\alpha\right)}^{-1}\left[\varphi_{\left(\alpha,\alpha\right)}\left(\Pi\left(\alpha\right)\right)+t\,\phi\left(\alpha\right)\right],\left(\alpha,\alpha\right)\right)\right),\label{expexp}
\end{equation}
being $\varphi_{\left(\alpha,\alpha\right)}^{-1}$ a local lateral
inverse of the immersion $\varphi_{\left(\alpha,\alpha\right)}$.
\end{prop}

It is natural to ask, given $\xi\in\mathfrak{g}$, if we can construct
$\exp\left(\xi\,t\right)$ up to quadratures. In the following subsection,
we shall give a partial answer to that question.

\subsubsection{The case of semisimple and compact Lie groups}

Let $G$ be a Lie group with Lie algebra $\mathfrak{g}.$

\begin{thm} \label{ibq} Consider $\xi\in\mathfrak{g}$ such that
\[
ad_{\xi}\left(\mathfrak{g}\right){}^{0}\cap\mathcal{R}_{Ad^{*}}\not=\emptyset,
\]
where $ad_{\xi}\left(\mathfrak{g}\right){}^{0}$ is the annihilator
in $\mathfrak{g}^{*}$ of the subspace $ad_{\xi}\left(\mathfrak{g}\right)\subseteq\mathfrak{g}$. 
\begin{enumerate}
\item If $\alpha_{0}\in ad_{\xi}\left(\mathfrak{g}\right){}^{0}\cap\mathcal{R}_{Ad^{*}}$,
then we can construct a Casimir $1$-form $\phi:\mathfrak{g}^{*}\to\mathfrak{g}$
such that $\phi\left(\alpha_{0}\right)=\xi.$ 
\item The curve $t\mapsto\exp\left(\xi\,t\right)$ can be obtained by quadratures. 
\end{enumerate}
\end{thm}

\noindent \textit{Proof. }Take $\alpha_{0}\in ad_{\xi}\left(\mathfrak{g}\right){}^{0}\cap\mathcal{R}_{Ad^{*}}$.
Then, $\alpha_{0}([\xi,\eta])=0,$ for all $\eta\in\mathfrak{g}$
or, in other words, $\xi\in\mathfrak{g}_{\alpha_{0}}$ with $\mathfrak{g}_{\alpha_{0}}$
the isotropy algebra of $\alpha_{0}$ with respect to the co-adjoint
representation of $G$ on $\mathfrak{g}^{*}$. Let $V\subseteq\mathfrak{g}^{*}$
be an admissible neighborhood of $\alpha_{0}$. Then (see Remark \ref{rem1})
the assigning 
\[
\alpha\in V\longmapsto\mathfrak{g}_{\alpha}\subseteq\mathfrak{g}
\]
defines a vector subbundle $W\coloneqq\coprod_{\alpha\in V}\mathfrak{g}_{\alpha}\rightarrow V$
of the trivial vector bundle $pr_{1}:V\times\mathfrak{g}\to V$. By
using the Inverse Function Theorem, we can construct an open subset
$\tilde{V}\subseteq V$ containing $\alpha_{0}$ and a section $\tilde{\phi}:\tilde{V}\subseteq\mathfrak{g}^{*}\to W$
(of such a bundle) satisfying $\tilde{\phi}\left(\alpha_{0}\right)=\xi.$
Note that, since $\tilde{\phi}\left(\alpha\right)\in\mathfrak{g}_{\alpha}$,
then $ad_{\tilde{\phi}\left(\alpha\right)}\alpha=0$, for all $\alpha\in\tilde{V}$.
Moreover, consider another open subsets $V_{1,2}$ such that $\alpha_{0}\in V_{1}\subseteq\bar{V}_{1}\subseteq V_{2}\subseteq\bar{V}_{2}\subseteq\tilde{V}$
and the bump function $\chi:\mathfrak{g}^{*}\rightarrow\mathbb{R}$
related to $V_{1,2}$, i.e. $\chi$ is equal to $1$ inside $\bar{V}_{1}$
and equal to $0$ outside $V_{2}$. It is clear that $\phi:\mathfrak{g}^{*}\to\mathfrak{g}$
given by 
\[
\phi\left(\alpha\right)=\left\{ \begin{array}{ll}
\tilde{\phi}\left(\alpha\right)\,\chi\left(\alpha\right), & \quad\alpha\in\tilde{V},\\
0, & \quad\alpha\notin\tilde{V},
\end{array}\right.
\]
satisfies the point $1.$ The point $2$ follows from $1$ and Proposition
\ref{p6} for $\alpha=\alpha_{0}$. $\;\;\;\blacksquare$\bigskip{}

For an important subclass of Lie groups, we have the following result.

\begin{thm}\label{quadrature-exponential} Let $G$ be a connected
Lie group with Lie algebra $\mathfrak{g}$ and $\mathcal{R}_{Ad}$
the open dense subset of $\mathfrak{g}$ which consists of the regular
points in $\mathfrak{g}$ with respect to the adjoint action of $G$
on $\mathfrak{g}.$ Suppose that there exists a non-degenerate ad-invariant
symmetric bilinear form $B:\mathfrak{g}\times\mathfrak{g}\to\mathbb{R}$.
Then, 
\begin{enumerate}
\item The linear map $B^{\flat}:\mathfrak{g}\to\mathfrak{g}^{*}$ given
by $\left\langle B^{\flat}\left(\xi\right),\eta\right\rangle =B\left(\xi,\eta\right)$,
for all $\xi,\eta\in\mathfrak{g}$, is a isomorphism satisfying $B^{\flat}\left(\mathcal{R}_{Ad}\right)=\mathcal{R}_{Ad^{*}}$,
and its inverse $B^{\sharp}:\mathfrak{g}^{*}\to\mathfrak{g}$ is a
Casimir $1$-form. 
\item For every $\xi\in\mathcal{R}_{Ad}$, the curve $t\mapsto\exp(\xi\,t)$
can be obtained by quadratures. 
\end{enumerate}
\end{thm}

\noindent \textit{Proof. } We have that (non-degeneracy) 
\begin{equation}
B\left(\xi,\eta\right)=0,\mbox{ for all }\eta\in\mathfrak{g}\quad\Longrightarrow\quad\xi=0,\label{non-degenerate}
\end{equation}
and (ad-invariance) 
\begin{equation}
B\left(\left[\xi,\eta\right],\nu\right)+B\left(\eta,\left[\xi,\nu\right]\right)=0,\mbox{ for all }\xi,\eta,\nu\in\mathfrak{g}.\label{ad-invariant}
\end{equation}
From (\ref{non-degenerate}), we deduce that $B^{\flat}$ is an isomorphism
of vector spaces. Moreover, using (\ref{ad-invariant}), it follows
that the following diagram

\begin{picture}(375,60)(60,40) \put(210,20){\makebox(0,0){${\mathfrak{g}}$}}\put(260,25){$B^{\flat}$}
\put(215,20){\vector(1,0){90}} \put(315,20){\makebox(0,0){${\mathfrak{g}}^{*}$}}
\put(190,50){$ad_{\xi}$} \put(210,30){\vector(0,1){40}}
\put(320,50){$ad_{\xi}^{*}$} \put(310,30){\vector(0,1){40}}
\put(210,80){\makebox(0,0){${\mathfrak{g}}$}} \put(260,85){$B^{\flat}$}
\put(215,80){\vector(1,0){90}} \put(315,80){\makebox(0,0){${\mathfrak{g}}^{*}$}}
\end{picture}

\vspace{1cm}

\noindent is commutative, for every $\xi\in{\mathfrak{g}}.$ So, since
$G$ is a connected Lie group, we also have that the diagram

\begin{picture}(375,60)(60,40) \put(210,20){\makebox(0,0){${\mathfrak{g}}$}}\put(260,25){$B^{\flat}$}
\put(215,20){\vector(1,0){90}} \put(315,20){\makebox(0,0){${\mathfrak{g}}^{*}$}}
\put(190,50){$Ad_{g}$} \put(210,30){\vector(0,1){40}} \put(320,50){$Ad_{g}^{*}=(Ad_{g^{-1}})^{*}$}
\put(310,30){\vector(0,1){40}} \put(210,80){\makebox(0,0){${\mathfrak{g}}$}}
\put(260,85){$B^{\flat}$} \put(215,80){\vector(1,0){90}}
\put(315,80){\makebox(0,0){${\mathfrak{g}}^{*}$}} \end{picture}

\vspace{1cm}

\noindent is commutative for every $g\in G.$ Thus, if $G_{\xi}$
(resp. $G_{B^{\flat}(\xi)}$) is the isotropy group of $\xi\in{\mathfrak{g}}$
(resp. $B^{\flat}(\xi)\in{\mathfrak{g}^{*}})$ with respect to the
adjoint (resp. co-adjoint) action of $G$ on ${\mathfrak{g}}$ (resp.
${\mathfrak{g}}^{*}$), we deduce that 
\[
G_{\xi}=G_{B^{\flat}\left(\xi\right)}.
\]
This implies that 
\begin{equation}
B^{\flat}\left(\mathcal{R}_{Ad}\right)=\mathcal{R}_{Ad^{*}}.\label{bads}
\end{equation}
On the other hand, from (\ref{ad-invariant}), we have that 
\[
B\left(\left[\xi,\eta\right],\xi\right)=0,\mbox{ for all }\xi,\eta\in\mathfrak{g}.
\]
Therefore, given $\alpha\in\mathfrak{g}^{*}$, if we write $\alpha=B^{\flat}\left(\xi\right)$
for some $\xi\in\mathfrak{g}$, we have for $B^{\sharp}=\left(B^{\flat}\right)^{-1}:\mathfrak{g}^{*}\to\mathfrak{g}$
that 
\[
\left\langle ad_{B^{\sharp}\left(\alpha\right)}^{*}\alpha,\eta\right\rangle =\left\langle \alpha,\left[B^{\sharp}\left(\alpha\right),\eta\right]\right\rangle =\left\langle B^{\flat}\left(\xi\right),\left[\xi,\eta\right]\right\rangle =B\left(\left[\xi,\eta\right],\xi\right)=0,
\]
for all $\eta\in\mathfrak{g}$. Then, $B^{\sharp}$ is a Casimir $1$-form.
This proves the first point. To prove the second point, note that,
according to \eqref{bads}, for every $\xi\in\mathcal{R}_{Ad}$ there
exists $\alpha\in\mathcal{R}_{Ad^{*}}$ such that $\xi=B^{\sharp}\left(\alpha\right)$.
Then, it is enough to use Proposition \ref{p6} for $\phi=B^{\sharp}$.
$\;\;\;\blacksquare$

\bigskip{}

\begin{rem*} It can be show that, under the conditions of the theorem
above, 
\[
ad_{\xi}(\mathfrak{g})^{0}\cap\mathcal{R}_{Ad^{*}}\not=\emptyset,\quad\forall\xi\in\mathcal{R}_{Ad}.
\]
So, the point $2$ of Theorem \ref{quadrature-exponential} can also
be proved by combining the equation above and Theorem \ref{ibq}.
\end{rem*}

Under the conditions of the last theorem, we can use \eqref{expexp}
for $\phi=B^{\sharp}$ and for all $\xi\in\mathcal{R}_{Ad}$, which
gives

\[
\exp\left(\xi\,t\right)=\pi_{G}\left(\left(\Pi,F\right)^{-1}\left(a_{\xi}\left(t\right),b_{\xi}\right)\right),
\]
with 
\[
\begin{array}{l}
a_{\xi}\left(t\right)=\varphi_{b_{\xi}}^{-1}\left[\varphi_{b_{\xi}}\left(\Pi\left(B^{\flat}\left(\xi\right)\right)\right)+t\,\xi\right]\quad\textrm{and}\quad b_{\xi}=F\left(B^{\flat}\left(\xi\right)\right)=\left(B^{\flat}\left(\xi\right),B^{\flat}\left(\xi\right)\right).\end{array}
\]

\begin{rem*} In particular, for $\xi\in\mathcal{R}_{Ad}$ and close
to $0$ (in order for $a_{\xi}\left(t\right)$ to be defined when
$t=1$), we have the following expression of the exponential map:
\[
\exp\left(\xi\right)=\pi_{G}\left(\left(\Pi,F\right)^{-1}\left(\varphi_{b_{\xi}}^{-1}\left[\varphi_{b_{\xi}}\left(\Pi\left(B^{\flat}\left(\xi\right)\right)\right)+\xi\right],b_{\xi}\right)\right).
\]

\end{rem*}

\begin{rem*} It is worth mentioning that $B^{\sharp}$ is an exact
$1$-form, i.e. $B^{\sharp}=dh$ with $h:\mathfrak{g}^{*}\rightarrow\mathbb{R}$
given by 
\[
h\left(\alpha\right)=\frac{1}{2}\,\left\langle \alpha,B^{\sharp}\left(\alpha\right)\right\rangle .
\]
Then, according to Remark \ref{ncir}, the related vector field $X_{\left|U\right.}^{B^{\sharp}}$
and the submersion $F$ define a NCI Hamiltonian system on $U$. \end{rem*}

For a semisimple Lie group $G$ with Lie algebra $\mathfrak{g}$,
the Killing form on $\mathfrak{g}$ satisfies the conditions in Theorem
\ref{quadrature-exponential} (see for example \cite{SW}). On the
other hand, a Lie algebra $\mathfrak{g}$ is the Lie algebra of a
compact Lie group if and only if $\mathfrak{g}$ admits an ad-invariant
scalar product (see, for instance, \cite{DK}). So, using Theorem
\ref{quadrature-exponential}, we have the next corollary.

\begin{cor} Let $G$ be a connected Lie group with Lie algebra $\mathfrak{g}$
and $\xi\in\mathcal{R}_{Ad}\subseteq\mathfrak{g}$. If $G$ is semisimple
or compact then $t\mapsto\exp(\xi\,t)$ can be obtained by quadratures.
\end{cor}

The last two results tell us that the exponential curve $\exp(\xi\,t)$
can be constructed by quadratures for $\xi$ living in an open dense
subset of $\mathfrak{g}$. Unfortunately, we can not ensure the same
for every Lie group.

\begin{rem} If $\mathfrak{g}$ is an arbitrary Lie algebra then the
subset 
\[
\left\{ \xi\in\mathfrak{g}\;:\;ad_{\xi}(\eta)^{0}\cap\mathcal{R}_{Ad^{*}}\not=\emptyset\right\} 
\]
is not, in general, dense in $\mathfrak{g}$. In fact, let ${\mathfrak{h}}(1,1)$
be the nilpotent Lie algebra of the Heisenberg group $H(1,1)$ of
dimension $3$. Then, we can consider a basis $\{\xi_{1},\xi_{2},\xi_{3}\}$
of ${\mathfrak{h}}(1,1)$ such that 
\[
[\xi_{1},\xi_{2}]=-[\xi_{2},\xi_{1}]=\xi_{3}
\]
and the rest of the basic Lie brackets are zero. So, if $(\alpha_{1},\alpha_{2},\alpha_{3})\in{\Bbb R}^{3}\cong{\mathfrak{h}}(1,1)^{*},$
we have that 
\[
\mathfrak{g}_{(\alpha_{1},\alpha_{2},\alpha_{3})}={\mathfrak{h}}(1,1),\:\:\mbox{ if \ensuremath{\alpha_{3}=0}},
\]
and 
\[
\mathfrak{g}_{(\alpha_{1},\alpha_{2},\alpha_{3})}=\left\langle \xi_{3}\right\rangle ,\;\;\mbox{ if }\alpha_{3}\not=0.
\]
Thus, we deduce that 
\[
\mathcal{R}_{Ad^{*}}=\left\{ {(\alpha_{1},\alpha_{2},\alpha_{3})}\in\mathbb{R}^{3}\;:\;\alpha_{3}\not=0\right\} ,
\]
which implies for $\xi=a^{1}\xi_{1}+a^{2}\xi_{2}+a^{3}\xi_{3}\in\mathfrak{g}$
that 
\[
ad_{\xi}\mathfrak{h}(1,1)^{0}\cap\mathcal{R}_{Ad^{*}}=\mathcal{R}_{Ad^{*}}\;\;\mbox{ if \ensuremath{a^{1}=a^{2}},}
\]
and 
\[
ad_{\xi}\mathfrak{h}(1,1)^{0}\cap\mathcal{R}_{Ad^{*}}=\emptyset\;\;\mbox{ if \ensuremath{a^{1}\not=a^{2}.}}
\]
\end{rem}

\subsection{Integrability conditions for invariant vertical fields}

Let us go back to Section \ref{ivf}. Consider a manifold $M$, a
vector field $X\in\mathfrak{X}\left(M\right)$ and a Lie group action
$\rho:G\times M\rightarrow M$. Assume that $X$ is vertical around
every point $m_{0}\in\mathcal{R}_{\rho}$ and $G$-invariant. Consider
a covering of $\mathcal{R}_{\rho}$ given by admissible neighborhoods
$U$, each one of them with an associated complete solution $\Sigma_{U}\coloneqq\rho\circ\left(id_{\Theta\left(U\right)}\times s\right)$,
as those given in Theorem \ref{1}, and the map 
\[
\eta_{U}:\pi\left(U\right)\rightarrow\mathfrak{g}\;:\;\lambda\mapsto\Theta_{*,s\left(\lambda\right)}\circ X\left(s\left(\lambda\right)\right)
\]
given by \eqref{el} in Theorem \ref{2}. From now on, we shall denote
$\mathfrak{g}_{\rho,m}$ the isotropy sub-algebra related to the point
$m$ and the action $\rho$.

\begin{thm}If for each $U$ and $\lambda\in\pi\left(U\right)$ we
have that 
\[
ad_{\eta_{U}(\lambda)+\varsigma_{\lambda}}(\mathfrak{g})^{0}\cap\mathcal{R}_{Ad^{*}}\not=\emptyset,
\]
for some $\varsigma_{\lambda}\in\mathfrak{g}_{\rho,s\left(\lambda\right)}$,
then $X$ is integrable up to quadratures along $\mathcal{R}_{\rho}$.
\end{thm}

\noindent \textit{Proof. } For a given $U$ and $\lambda\in\pi(U)$,
we know that the integral curves of $X$, with initial conditions
inside $U$, are given by the formula (see \eqref{Gtrg}) 
\[
\Gamma\left(t\right)=\rho\left(g_{0}\,\exp\left(\left(\eta_{U}\left(\lambda\right)+\chi_{\lambda}\right)\,t\right),s\left(\lambda\right)\right),
\]
with $\lambda\in\pi(U)$ and $\chi_{\lambda}\in\mathfrak{g}_{\rho,s\left(\lambda\right)}$
arbitrary. On the other hand, using Theorem \ref{ibq}, given 
\[
\alpha\in ad_{\eta_{U}(\lambda)+\varsigma_{\lambda}}(\mathfrak{g})^{0}\cap\mathcal{R}_{Ad^{*}},
\]
we can construct a Casimir $1$-form $\phi_{\lambda}:\mathfrak{g}^{*}\to\mathfrak{g}$
such that $\phi_{\lambda}(\alpha)=\eta_{U}(\lambda)+\varsigma_{\lambda}$.
Thus, taking $\chi_{\lambda}=\varsigma_{\lambda}$, we have that 
\[
\Gamma(t)=\rho(g_{0}\exp(\phi_{\lambda}(\alpha)\,t),s(\lambda))
\]
and, using Proposition \ref{p6} for $\phi=\phi_{\lambda}$, it follows
that $\Gamma$ can be constructed up to quadratures. $\;\;\;\blacksquare$

\bigskip{}

Now, let us suppose that $M$ is a symplectic manifold, with symplectic
structure $\omega$, and $\rho:G\times M\to M$ is a symplectic action
with $Ad^{*}$-equivariant momentum map $K:M\rightarrow\mathfrak{g}^{*}$.

\begin{thm} \label{kar} Consider an equivariant function $\phi:\mathfrak{g}^{*}\rightarrow\mathfrak{g}$
and the vector field $X=\omega^{\sharp}\circ K^{*}\phi$. If 
\[
K\left(\mathcal{R}_{\rho}\right)\cap\mathcal{R}_{Ad^{*}}\neq\emptyset
\]
and $\phi$ is also a Casimir $1$-form, then there exists a $G$-invariant
open subset of $V\subseteq\mathcal{R}_{\rho}$ where $X$ is integrable
up to quadratures. \end{thm}

\noindent \textit{Proof. } Under above condition, according to Proposition
\ref{7}, $X=\omega^{\sharp}\circ K^{*}\phi$ is vertical and $G$-invariant.
On the other hand, if $K\left(\mathcal{R}_{\rho}\right)\cap\mathcal{R}_{Ad^{*}}\neq\emptyset$,
according to Proposition \ref{5}, there exists a $G$-invariant open
subset $V\subseteq\mathcal{R}_{\rho}$ such that $K\left(V\right)\subseteq\mathcal{R}_{Ad^{*}}$.
Consider a covering of $V$ as above and the related maps $\eta_{U}$.
Using Proposition \ref{7} again, 
\[
\eta_{U}\left(\lambda\right)=\phi\left(K\left(s\left(\lambda\right)\right)\right)+\xi_{\lambda},
\]
for some $\xi_{\lambda}\in\mathfrak{g}_{\rho,s\left(\lambda\right)}$,
and the integral curves of $X$ by points of $U$ are of the form
\begin{equation}
\Gamma(t)=\rho\left(g\exp\left(\phi\left(K\left(s\left(\lambda\right)\right)\right)\,t\right),s\left(\lambda\right)\right),\;\;\mbox{ with \ensuremath{\lambda\in\Pi\left(U\right)} and \ensuremath{g\in\Theta(U)}}.\label{5.27'}
\end{equation}
Then, since $s\left(\lambda\right)\in U\subseteq V$, it follows that
$K\left(s\left(\lambda\right)\right)\in\mathcal{R}_{Ad^{*}}$, and
consequently, using (\ref{5.27'}) and Proposition \ref{p6} for $\alpha=K\left(s\left(\lambda\right)\right)$,
we deduce the result. $\;\;\;\blacksquare$

\bigskip{}
 More interesting examples can be constructed by using the next lemma.

\begin{lem} If $h:\mathfrak{g}^{*}\rightarrow\mathbb{R}$ is a $G$-invariant
function with respect to $Ad^{*}$, then $dh:\mathfrak{g}^{*}\rightarrow\mathfrak{g}$
is equivariant and a Casimir $1$-form. \end{lem}

For a proof, see \cite{adler}, Lemma 2.9.

\begin{thm} Consider $G$-invariant functions $h_{i}:\mathfrak{g}^{*}\rightarrow\mathbb{R}$
(resp. $f_{i}:M\rightarrow\mathbb{R}$), $i=1,...,k$, with respect
to $Ad^{*}$ (resp. $\rho$). Suppose that $K\left(\mathcal{R}_{\rho}\right)\cap\mathcal{R}_{Ad^{*}}\neq\emptyset$
and define 
\begin{equation}
X\left(m\right)=\sum_{i=1}^{k}f_{i}\left(m\right)\,\left(\omega^{\sharp}\circ K^{*}dh_{i}\right)\left(m\right),\;\;\;\forall m\in M.\label{xfihi}
\end{equation}
Then, there exists a $G$-invariant open subset $V\subseteq\mathcal{R}_{\rho}$
where the vector field $X$ is integrable up to quadratures. \end{thm}

\noindent \textit{ Proof. }Since each field $\omega^{\sharp}\circ K^{*}dh_{i}$
is $G$-invariant and vertical, the same is true for $X$. On the
other hand, given (as in the proof of Theorem \eqref{kar}) a $G$-invariant
open subset of $V\subseteq\mathcal{R}_{\rho}$ such that $K\left(V\right)\subseteq\mathcal{R}_{Ad^{*}}$,
a covering of $V$ by admissible neighborhoods $U$ and the related
maps $\eta_{U}$, for each $\lambda\in\pi\left(U\right)$ we have
that 
\[
\eta_{U}\left(\lambda\right)=\sum_{i=1}^{k}f_{i}\left(s\left(\lambda\right)\right)\,dh_{i}\left(K\left(s\left(\lambda\right)\right)\right)+\xi_{\lambda},
\]
for some $\xi_{\lambda}\in\mathfrak{g}_{\rho,s\left(\lambda\right)}$.
Then, defining $\phi_{\lambda}:\mathfrak{g}^{*}\rightarrow\mathfrak{g}$
by 
\[
\phi_{\lambda}\left(\alpha\right)=\sum_{i=1}^{k}f_{i}\left(s\left(\lambda\right)\right)\,dh_{i}\left(\alpha\right),
\]
which is a Casimir $1$-form, we have that 
\[
\eta_{U}\left(\lambda\right)=\phi_{\lambda}\left(K\left(s\left(\lambda\right)\right)\right)+\xi_{\lambda}.
\]
Finally, since $K\left(s\left(\lambda\right)\right)\in\mathcal{R}_{Ad^{*}}$
(as we saw in the previous theorem), the theorem follows from Proposition
\ref{p6} for $\phi=\phi_{\lambda}$ and $\alpha=K\left(s\left(\lambda\right)\right)$.
$\;\;\;\blacksquare$

\bigskip{}

It is worth mentioning that the vector field $X$ given by \eqref{xfihi}
is not, in general, a Hamiltonian vector field.

\section*{Acknowledgements}

This work was partially supported by Ministerio de Ciencia Innovación
(Spain), FEDER co-financing, grants PGC2018-098265-B-C32 (JC M. and
E.P.) and CONICET (S.G.).

\end{document}